\documentclass[11pt,a4paper]{article}
\pdfoutput=1 % for arxiv

\usepackage{amsmath,amssymb,amsthm,mathtools}
\usepackage{color,url,comment}
\usepackage[T1]{fontenc}
\usepackage{lmodern}
\usepackage{enumitem}
\usepackage{cite}
\usepackage[colorlinks=true,linkcolor=blue,citecolor=blue,urlcolor=blue]{hyperref}

\allowdisplaybreaks[1]

\newtheorem{definition}{Definition} [section]
\newtheorem{theorem}[definition]{Theorem}
\newtheorem{lemma}[definition]{Lemma}
\newtheorem{proposition}[definition]{Proposition}

\newtheorem{corollary}[definition]{Corollary}
\newtheorem{conjecture}[definition]{Conjecture}
\newtheorem{claim}[definition]{Claim}

\newtheorem{construction}[definition]{Construction}

\theoremstyle{remark}

\newenvironment{poc}{\begin{proof}[Proof of the claim]}{\end{proof}}
\def\ex{\mathrm{ex}}
\def\coex{\mathrm{coex}}

\newcommand{\hide}[1]{}

\usepackage{tikz}
\usetikzlibrary{arrows.meta,calc}
\usepackage{caption}

\definecolor{Ink}{HTML}{243443}
\definecolor{Muted}{HTML}{657580}
\definecolor{Border}{HTML}{DCE3E8}
\definecolor{Label}{HTML}{247DA5}
\definecolor{Test}{HTML}{CE752B}
\definecolor{TestAlt}{HTML}{B8493B}
\definecolor{Compare}{HTML}{0072B2}
\definecolor{CompareAlt}{HTML}{CC3366}
\definecolor{Reservoir}{HTML}{586CB0}
\definecolor{Filler}{HTML}{188577}
\definecolor{Cycle}{HTML}{526C7F}
\pagecolor{white}
\color{Ink}
\tikzset{
	part/.style={circle,fill=white,draw=Ink,line width=1.15pt,
		minimum size=10mm,inner sep=1pt,font=\large,text=Ink},
	labelled/.style={part,draw=Label},
	testpart/.style={part,draw=Test},
	comparepart/.style={part,draw=Compare},
	comparecpart/.style={part,draw=CompareAlt},
	reservoirpart/.style={part,draw=Reservoir},
	fillerpart/.style={part,draw=Filler},
	heading/.style={font=\sffamily\bfseries\normalsize,text=Ink}
}

\begin{document}
	\title{\bf\Large The spectrum of codegree Turán densities of 3-graphs}
	\date{}
	\author{Jun Gao\thanks{Academy of Mathematics and Systems Science, Chinese Academy of Sciences, Beijing 100084, China. Email: \texttt{jungao@amss.ac.cn}}
		\and Oleg Pikhurko\thanks{Mathematics Institute and DIMAP, University of Warwick, Coventry, UK. Research supported by ERC Advanced Grant 101020255 and a grant of access to OpenAI models through the ChatGPT for Academic Researchers program. Email: \texttt{o.pikhurko@warwick.ac.uk}}
		\and Mingyuan Rong\thanks{School of Mathematical Sciences, USTC, Hefei, China. Supported by National Key R\&D Program of China 2023YFA1010201, NSFC Grant No. 12125106 and the USTC Excellent PhD Students Overseas Study Program. Email: \texttt{rong\_ming\_yuan@mail.ustc.edu.cn}}
		\and Shumin Sun\thanks{School of Mathematics, University of Birmingham, Birmingham, UK. Supported by the EPSRC New Investigator Award EP/Y006399/1. Email: \texttt{s.sun.5@bham.ac.uk}}}
	% Required for inserting images
	
	\maketitle

	\begin{abstract}
		The \emph{minimum codegree} $\delta_{k-1}(H)$ of a $k$-uniform hypergraph (in short, $k$-graph) $H$ is the smallest number of edges containing a $(k-1)$-set of vertices. For a family $\mathcal{F}$ of $k$-graphs, its \emph{codegree Turán density} $\gamma(\mathcal{F})$ is the asymptotic maximum of $\delta_{k-1}(H)/|V(H)|$ over growing $\mathcal{F}$-free $k$-graphs $H$. 
		
		We prove that every real number $\alpha\in[0,1)$ is the codegree Turán density of an infinite family $\mathcal{F}_\alpha$ of $3$-graphs, and this family can be chosen so that $\gamma(\mathcal{F}')>\alpha$ for every finite subfamily $\mathcal{F}'\subseteq\mathcal{F}_\alpha$. This resolves the $3$-graph case of a conjecture of Mubayi and Zhao, the only case that remained open after recent work of Lin and Zhou. (Li and Liu earlier determined the codegree Turán spectrum in every uniformity $k\ge 4$ and, very recently, they determined the spectrum for $k=3$  independently of us.)

		We also consider \emph{positive codegree Turán density}, defined analogously by taking the minimum only over $(k-1)$-sets contained in at least one edge. We prove that there is an absolute constant $c<1$ such that, for every integer $k\geq3$, every real number in $[c,1)$ is the positive codegree Tur\'an density of a $k$-graph family. This is done by presenting a construction with a relatively short proof, without optimising the value of~$c$.
	\end{abstract}
	
	\section{Introduction}
	
	For an integer $k\geq 2$, a \emph{$k$-uniform hypergraph}, or simply a \emph{$k$-graph}, is a pair $H=(V(H),E(H))$ with
	$E(H)\subseteq\binom{V(H)}{k}$. 
	All hypergraphs considered in this paper are finite and simple. 
	Given a possibly infinite family $\mathcal{F}$ of
	$k$-graphs, we say that $H$ is \emph{$\mathcal{F}$-free} if it contains no member of $\mathcal{F}$ as a subgraph.
	The \emph{Tur\'an number} $\ex(n,\mathcal{F})$ is the
	largest number of edges in an $n$-vertex $\mathcal{F}$-free $k$-graph, and the corresponding \emph{Tur\'an density} is
	\[
	\pi(\mathcal{F})
	\coloneqq \lim_{n\to\infty}
	\frac{\ex(n,\mathcal{F})}{\binom{n}{k}}.
	\]
	
	The existence of this limit follows from a standard averaging argument.
	Determining Tur\'an densities is a central problem in extremal
	combinatorics, and even basic examples remain difficult in uniformity at least three. 
	We refer the reader to the surveys~\cite{kee11,sido95} for
	background.
	If $\mathcal F$ contains an edgeless hypergraph, then
	$\ex(n,\mathcal F)$ and its various versions are undefined for all sufficiently large $n$.
	We therefore assume throughout that every member of a forbidden family contains at least one edge.
	
	Our main focus is the codegree Tur\'an problem, introduced by Mubayi and Zhao~\cite{MZ07}. For a $k$-graph $H$ and a set
	$S\in\binom{V(H)}{k-1}$, write
	\[
	N_H(S)\coloneqq
	\{v\in V(H)\setminus S:S\cup\{v\}\in E(H)\},
	\qquad \deg_H(S)\coloneqq |N_H(S)|.
	\]
	The \emph{minimum codegree} of $H$ is
	\[
	\delta_{k-1}(H)
	\coloneqq \min_{S\in\binom{V(H)}{k-1}} \deg_H(S).
	\]
	For a family $\mathcal{F}$ of $k$-graphs, define its \emph{codegree Tur\'an number} by
	\[
	\coex(n,\mathcal{F})
	\coloneqq
	\max\{\delta_{k-1}(H): |V(H)|=n
	\text{ and }H\text{ is }\mathcal{F}\text{-free}\},
	\]
	and its \emph{codegree Tur\'an density} by
	\[
	\gamma(\mathcal{F})
	\coloneqq \lim_{n\to\infty}
	\frac{\coex(n,\mathcal{F})}{n}.
	\]
	Mubayi and Zhao~\cite[Proposition~1.2]{MZ07} proved that this limit exists for every family $\mathcal{F}$, including infinite families. When $\mathcal{F}=\{F\}$, we write $\gamma(F)$ instead of $\gamma(\{F\})$, and similarly for the other parameters.
	
	A basic question is which real numbers occur as Tur\'an densities. For each $k\geq 2$, let
	\[
	\Pi_k\coloneqq
	\{\pi(\mathcal{F}):\mathcal{F}
	\text{ is a nonempty family of }k\text{-graphs}\}.
	\]
	The Erd\H{o}s--Stone--Simonovits theorem~\cite{ES46,ES66} gives
	\[
	\Pi_2=\left\{1-\frac{1}{r}:r=1,2,\ldots\right\}.
	\]
	In particular, every point of $[0,1)$ has a nonempty open interval on its right containing
	no graph Tur\'an density. More generally, a number $a\in[0,1)$ is called
	a \emph{jump} for $k$-graphs if some $\varepsilon>0$ satisfies
	$\Pi_k\cap(a,a+\varepsilon)=\varnothing$. The existence of non-jumps in
	higher uniformities, established by Frankl and R\"odl~\cite{FR84}, shows
	that the structure of hypergraph Tur\'an densities differs substantially
	from the graph case. For further examples of non-jumps, see e.g.~\cite{FPRT07,P14,YP23,LM26,Shaw26}. On the other hand, Erd\H{o}s~\cite{E64} showed that, for every
	$k\geq 3$, all numbers in $[0,k!/k^k)$ are jumps for $k$-graphs. For $k=3$, Baber and Talbot~\cite{BT11} later found further
	jumps beyond this interval.
	% See also~\cite{P14} for further developments	on these spectra. 
	Recently, Liu and Pikhurko~\cite{LP26} showed that, for every
	$k\geq 3$, $\Pi_k$ contains an interval $[c_k,1)$ for some $c_k<1$, and Li and Liu~\cite{LiLiu26v2} proved that $\Pi_3\supseteq [4/9,1)$. 
	
	%\op{Add other references: BaberTalbot10 on further jumps, LiuPikhurko26 on interval of Turan densities, MubayLiu that 4/9 is a non-jump}
	
	We consider the corresponding question for codegree densities. Define
	\[
	\Gamma_k\coloneqq
	\{\gamma(\mathcal{F}):\mathcal{F}
	\text{ is a nonempty family of }k\text{-graphs}\}
	\subseteq[0,1).
	\]
	For graphs, the degree form of the Erd\H{o}s--Stone--Simonovits theorem implies that $\Gamma_2=\Pi_2$. For every $k\geq 3$, however, Mubayi and Zhao~\cite[Theorem~1.6]{MZ07} proved that $\Gamma_k$ is dense in $[0,1)$. 
	They asked whether $\Gamma_k=[0,1)$ and proposed the following stronger conjecture.
	
	\begin{conjecture}[Mubayi--Zhao~\cite{MZ07}]
		\label{conj:MZ}
		For every integer $k\geq 3$ and every real number $\alpha\in[0,1)$,
		there exists an infinite family $\mathcal{F}$ of $k$-graphs such that
		\[
		\gamma(\mathcal{F})=\alpha,
		\qquad\text{and}\qquad
		\gamma(\mathcal{F}')>\alpha
		\quad\text{for every finite }\mathcal{F}'\subseteq\mathcal{F}.
		\]
	\end{conjecture}
	
	Since there are only countably many finite families of $k$-graphs up to isomorphism, realizing every real number in $[0,1)$ requires allowing infinite forbidden families. On the other hand, in our earlier work~\cite{GPRS26}, we proved that every rational number in $[0,1)$ is the codegree Tur\'an density of a finite family of $k$-graphs, for every $k\geq3$.
	
	Recently, Li and Liu~\cite{LiLiu2026} determined the minimum $\ell$-degree Tur\'an density spectrum for all $k>\ell\geq 2$ except $(k,\ell)=(3,2)$, that is, the case of codegree
	Tur\'an density for $3$-graphs. Lin and Zhou~\cite{LinZhou2026} established this together with the stronger finite-subfamily conclusion for all $k>\ell>2$, thus proving Conjecture~\ref{conj:MZ} for every $k\geq 4$.
	
	Our main result (Theorem~\ref{thm:spectrum} below) proves Conjecture~\ref{conj:MZ} for $k=3$, the only remaining open case. When we were in the final stages of writing up this paper, Version 2 of Li and Liu~\cite{LiLiu26v2} appeared, also proving that $\Gamma_3=[0,1)$.
	
	\begin{theorem}
		\label{thm:spectrum}
		For every real number $\alpha\in[0,1)$, there exists an infinite family
		$\mathcal{F}_\alpha$ of $3$-graphs such that
		\[
		\gamma(\mathcal{F}_\alpha)=\alpha.
		\]
		Moreover, $\mathcal{F}_\alpha$ can be chosen so that
		\[
		\gamma(\mathcal{F}')>\alpha
		\qquad\text{for every finite }
		\mathcal{F}'\subseteq\mathcal{F}_\alpha.
		\]
	\end{theorem}
	
	We also give an explicit operation that sends a family of $k$-graphs to a family of $(k+1)$-graphs with the same codegree Tur\'an density. The operation also preserves the property that every finite subfamily has strictly larger density.
	
	\begin{theorem}\label{thm:lift-uniformity}
		Let $k\geq2$ be an integer, and let $\mathcal F$ be a family of $k$-graphs with $\gamma(\mathcal F)=\alpha$ such that
		\[
		\gamma(\mathcal F_0)>\alpha
		\qquad\text{for every finite }
		\mathcal F_0\subseteq\mathcal F.
		\]
		Then there exists a family $\mathcal B$ of $(k+1)$-graphs with
		$\gamma(\mathcal B)=\alpha$ such that
		\[
		\gamma(\mathcal B_0)>\alpha
		\qquad\text{for every finite }
		\mathcal B_0\subseteq\mathcal B.
		\]
	\end{theorem}
	
	Together with Theorem~\ref{thm:spectrum}, this theorem yields
	another proof of Conjecture~\ref{conj:MZ} for $k\geq 4$.
	
	A related question concerns positive codegree Tur\'an densities, where the minimum is taken only over sets contained in an edge.
	More specifically, for a nonempty $k$-graph $H$, define its \emph{minimum positive codegree} by
	\[
	\delta^+_{k-1}(H)
	:=\min\{\deg_H(S):S\in\tbinom{V(H)}{k-1},\ \deg_H(S)>0\},
	\]
	and set $\delta^+_{k-1}(H)=0$ when $H$ is edgeless. For a nonempty,
	possibly infinite family $\mathcal F$ of $k$-graphs, each having an edge,
	put
	\[
	\begin{aligned}
		\operatorname{ex}^+(n,\mathcal F)
		&:=\max\{\delta^+_{k-1}(H):|V(H)|=n,\ H\text{ is }\mathcal F\text{-free}\},\\
		\gamma^+(\mathcal F)&:=\lim_{n\to\infty}
		\frac{\operatorname{ex}^+(n,\mathcal F)}n,
	\end{aligned}
	\]
	where the limit exists by~\cite[Proposition~1]{Pikhurko23} (see also~\cite[Theorem~20]{HLP25} for the case of a single forbidden $k$-graph).
	
	%\op{HalfpapLemonsPalmer25 also proved this} 
	
	Our second main result shows that every number sufficiently close to one from below occurs as a positive codegree Tur\'an density, in every uniformity at least three.
	\begin{theorem}\label{thm:positive-introduction}
		There is an absolute constant $c<1$ such that, for every integer $k\ge3$ and every real number
		$\alpha\in[c,1)$, there is a nonempty, possibly infinite family
		$\mathcal F$ of $k$-graphs, each member having an edge, such that
		$\gamma^+(\mathcal F)=\alpha$.
	\end{theorem}
	
	Note that one cannot expect that positive codegree densities of $k$-graphs span the whole interval $[0,1]$. Indeed, Halfpap, Lemons and Palmer~\cite[Corollary 25]{HLP25} showed that no $k$-graph $F$ can satisfy $0<\gamma^+(F)<1/k$, and their arguments can be easily extended to arbitrary families of $k$-graphs. Furthermore, Balogh, Halfpap, Lidick\'y and Palmer~\cite[Theorem 2.1]{BaloghHalfpapLidickyPalmer26} proved that the interval $(1/k,2/(2k-1))$ is also free of positive codegree densities of $k$-graphs.
	
	The rest of the paper is organised as follows.
	In Section~\ref{sec:interval}, we construct families of $3$-graphs
	whose codegree Tur\'an densities cover an interval, and prove
	the matching lower and upper bounds.
	In Section~\ref{sec:spectrum}, we extend this interval to $(0,1)$
	via a cyclic substitution based on the construction
	in~\cite[Proof of Claim~2.3]{GPRS26}.
	In Section~\ref{sec:finite-gaps}, we establish the finite-subfamily
	conclusion through a further refinement of the construction
	and a separate treatment of density zero. In Section~\ref{sec:lift-uniformity}, we introduce the operation that can lift a family of $k$-graphs to a family of $(k+1)$-graphs with the same codegree Tur\'an density, and prove Theorem~\ref{thm:lift-uniformity}.
	Finally, in Section~\ref{sec:positive-terminal}, we prove Theorem~\ref{thm:positive-introduction} by constructing a short interval of positive codegrees via a 3-graph core construction and amplifying it directly in every uniformity. The proof gives $c=1-72/(2^{72}+152)$; our aim here is simplicity rather than a small left endpoint.
	\medskip
	
	We will use the following notation. We will write $[n]$ for the set $\{1,\dots,n\}$. 
	%For a $k$-graph $G$, $v(G):=|V(G)|$ denotes its \emph{order}. 
	A set $I\subseteq V(G)$ is \emph{independent} in a $k$-graph $G$ if no edge of $G$ lies entirely inside $I$. We let $K_m^{(k)}$ denote the complete $k$-graph with $m$ vertices.

	\section{An interval of codegree Tur\'an densities}\label{sec:interval}
	In this section, we first obtain explicit bounds for families realizing an interval of small codegree densities.
	
	\begin{theorem}
		\label{thm:initial-interval}
		For every integer $s\geq 0$ and every $\theta\in[0,1]$, there exists a
		family $\mathcal{F}_{s,\theta}$ of $3$-graphs such that, for
		every integer $n\geq 9+s$,
		\[
		\frac{n}{8+s-\theta}-3
		\leq \coex(n,\mathcal{F}_{s,\theta})
		\leq \frac{n}{8+s-\theta}.
		\]
		In particular,
		\[
		\gamma(\mathcal{F}_{s,\theta})
		=\frac{1}{8+s-\theta}.
		\]
	\end{theorem}
	
	The densities in Theorem~\ref{thm:initial-interval} cover $(0,1/7]$.
	Since $0$ is the codegree density when we forbid a single edge, this gives the following corollary.
	
	\begin{corollary}\label{cor:base}
		For every $\beta\in[0,1/7]$, there exists a family
		$\mathcal F_\beta$ of $3$-graphs such that
		\[
		\gamma(\mathcal F_\beta)=\beta.
		\]
	\end{corollary}
	
	\subsection{Subgraph-closed family and the construction}\label{sec:base}
	
	A family $\mathcal C$ of $3$-graphs is \emph{subgraph-closed} if it is closed under deleting vertices and edges. For such a class containing every edgeless $3$-graph, define
	\[
	c_{\mathcal C}(n)
	\coloneqq\max\{\delta_2(G):G\in\mathcal C,\ |V(G)|=n\}.
	\]
	Let $\mathcal F(\mathcal C)$ be the family of all finite $3$-graphs
	not in $\mathcal C$.
	Since $\mathcal C$ is closed under taking subgraphs,
	a $3$-graph is $\mathcal F(\mathcal C)$-free if and only if
	it belongs to $\mathcal C$. Hence
	\begin{equation}\label{eq:class}
		\coex(n,\mathcal F(\mathcal C))=c_{\mathcal C}(n).
	\end{equation}
	We define
	\[
	\operatorname{dens}(\mathcal C)
	\coloneqq\lim_{n\to\infty}\frac{c_{\mathcal C}(n)}n
	= \lim_{n\to\infty}\frac{\coex(n,\mathcal F(\mathcal C))}n= \gamma(\mathcal F(\mathcal C)),
	\]
	where the existence of the limit, as we discussed earlier, follows from~\cite[Proposition~1.2]{MZ07}.
	
	We identify the circle with $[0,1)$, with addition modulo one. For
	$x\in[0,1)$ and $\theta\in[0,1]$, define
	\[
	I_\theta(x)
	\coloneqq\{y\in[0,1):0<(y-x)\bmod1<\theta\}.
	\]
	In other words, $I_\theta(x)$ is the open arc of length $\theta$
	immediately following $x$ in the positive direction. In particular, $I_0(x)=\varnothing$ and $I_1(x)=[0,1)\setminus\{x\}$.

Let us present the core of our construction, based on the one in Mubayi and Zhao~\cite[Section~3.2]{MZ07} (that basically uses only Rules 1 and 2). 

	\begin{construction}[$\mathcal A_{s,\theta}$ and $\mathcal F_{s,\theta}$]\label{cons:s-theta}
		Fix an integer $s\geq0$ and $\theta\in[0,1]$, and take pairwise
		disjoint vertex sets
		\[
		X,Y,Z,A,B,C,D,W_1,\ldots,W_s,R.
		\]
		We call the sets other than $R$ \emph{basic sets} and write
		$m\coloneqq7+s$ for their number.
		Assign each vertex $v\in X\cup Y\cup Z$ a label
		$\lambda(v)\in[0,1)$; we call these sets (and their vertices) \emph{labelled}.
		Empty sets and repeated labels are allowed.
		
		Arrange the basic sets in the cyclic order
		\begin{equation}\label{eq:basic-cycle}
			(V_0,\ldots,V_{m-1})
			\coloneqq(X,Y,Z,A,B,C,D,W_1,\ldots,W_s).
		\end{equation}
		The six pairs $\{U,B\}$ and $\{U,C\}$ with $U\in\{X,Y,Z\}$
		are called \emph{comparison pairs}, and $\{X,A\}$ and $\{Y,D\}$
		are called \emph{test pairs}.
		These eight pairs are \emph{special}, and all remaining pairs
		of distinct sets are \emph{ordinary}.
		
		For each choice of vertex sets and labels, define a host $3$-graph
		by the following edge rules, with no other triples included.
		The notation $V_iV_iV_{i+1}$ denotes a triple with two vertices
		in $V_i$ and one in $V_{i+1}$, with indices taken modulo $m$.
		\begin{enumerate}
			\item Include every triple of type $V_iV_iV_{i+1}$ and every
			triple of type $RRB$.
			
			\item Include every triple meeting three distinct sets if all three
			pairs of its sets are ordinary.
			
			\item For two vertices in distinct labelled sets and a third
			vertex in $B$ or $C$, include their triple if and only if the
			corresponding inequality in the following table holds.
			\begin{equation}\label{eq:comparison-table}
				\begin{array}{c|cc}
					\text{labelled vertices}
					&\text{third vertex in }B
					&\text{third vertex in }C\\\hline
					x\in X,\ y\in Y
					&\lambda(y)<\lambda(x)&\lambda(x)<\lambda(y)\\
					y\in Y,\ z\in Z
					&\lambda(z)<\lambda(y)&\lambda(y)<\lambda(z)\\
					z\in Z,\ x\in X
					&\lambda(x)<\lambda(z)&\lambda(z)<\lambda(x)
				\end{array}
			\end{equation}
			
			\item Include every triple of type $XAR$ or $YDR$.
			For $x\in X$ and $y\in Y$, also include $xya$ with $a\in A$
			if $\lambda(y)\in I_\theta(\lambda(x))$, and include $xyd$
			with $d\in D$ if $\lambda(x)\in I_\theta(\lambda(y))$.
		\end{enumerate}
		
		Let $\mathcal A_{s,\theta}$ consist of all subgraphs of these hosts,
		allowing all set sizes and labels.
		Define
		\[
		\mathcal F_{s,\theta}
		\coloneqq\mathcal F(\mathcal A_{s,\theta}).
		\]
	\end{construction}

	\begin{center}
		\begin{figure}
			\resizebox{\linewidth}{!}{%
				\begin{tikzpicture}[x=1cm,y=1cm,font=\normalsize,line cap=round,line join=round]
					% Oval radii: 2.85 cm horizontally, 3.60 cm vertically.
					% The fixed 17.2 cm width preserves all text, node and stroke sizes.
					\path[use as bounding box] (-.10,0) rectangle (17.10,11.3630);
					% Version P5: symbolic W-chain, no oval guide contour.
					%\node[heading] at (4.6500,10.4100) {Cyclic types and special pairs};
					\coordinate (p5-X) at (2.3443,7.7975);
					\coordinate (p5-Y) at (3.7693,9.1053);
					\coordinate (p5-Z) at (5.5307,9.1053);
					\coordinate (p5-A) at (6.9557,7.7975);
					\coordinate (p5-B) at (7.5000,5.6815);
					\coordinate (p5-C) at (6.9557,3.5655);
					\coordinate (p5-D) at (5.5307,2.2577);
					\coordinate (p5-Wone) at (3.7693,2.2577);
					\coordinate (p5-Wlast) at (1.8000,5.6815);
					\coordinate (p5-R) at (9.0500,5.6815);
					\draw[draw=Compare,line width=1.1pt] (p5-X)--(p5-B);
					\draw[draw=CompareAlt,line width=1.1pt] (p5-X)--(p5-C);
					\draw[draw=Compare,line width=1.1pt] (p5-Y)--(p5-B);
					\draw[draw=CompareAlt,line width=1.1pt] (p5-Y)--(p5-C);
					\draw[draw=Compare,line width=1.1pt] (p5-Z)--(p5-B);
					\draw[draw=CompareAlt,line width=1.1pt] (p5-Z)--(p5-C);
					\draw[draw=Test,line width=1.25pt,dash pattern=on 4pt off 2.4pt] (p5-X)--(p5-A);
					\draw[draw=Test,line width=1.25pt,dash pattern=on 4pt off 2.4pt] (p5-Y)--(p5-D);
					\draw[draw=Cycle,line width=1.05pt,-{Latex[length=2.8mm,width=1.9mm]},shorten <=5.35mm,shorten >=5.35mm] (p5-X) arc[start angle=144,end angle=108,x radius=2.85cm,y radius=3.6000cm];
					\node[anchor=east,font=\small,text=Cycle,fill=white,inner sep=1.2pt] at (2.7632,8.7956) {$XXY$};
					\draw[draw=Cycle,line width=1.05pt,-{Latex[length=2.8mm,width=1.9mm]},shorten <=5.35mm,shorten >=5.35mm] (p5-Y) arc[start angle=108,end angle=72,x radius=2.85cm,y radius=3.6000cm];
					\node[anchor=south,font=\small,text=Cycle,fill=white,inner sep=1.2pt] at (4.6500,9.5307) {$YYZ$};
					\draw[draw=Cycle,line width=1.05pt,-{Latex[length=2.8mm,width=1.9mm]},shorten <=5.35mm,shorten >=5.35mm] (p5-Z) arc[start angle=72,end angle=36,x radius=2.85cm,y radius=3.6000cm];
					\node[anchor=west,font=\small,text=Cycle,fill=white,inner sep=1.2pt] at (6.5368,8.7956) {$ZZA$};
					\draw[draw=Cycle,line width=1.05pt,-{Latex[length=2.8mm,width=1.9mm]},shorten <=5.35mm,shorten >=5.35mm] (p5-A) arc[start angle=36,end angle=0,x radius=2.85cm,y radius=3.6000cm];
					\node[anchor=west,font=\small,text=Cycle,fill=white,inner sep=1.2pt] at (7.7029,6.8710) {$AAB$};
					\draw[draw=Cycle,line width=1.05pt,-{Latex[length=2.8mm,width=1.9mm]},shorten <=5.35mm,shorten >=5.35mm] (p5-B) arc[start angle=0,end angle=-36,x radius=2.85cm,y radius=3.6000cm];
					\node[anchor=west,font=\small,text=Cycle,fill=white,inner sep=1.2pt] at (7.7029,4.4920) {$BBC$};
					\draw[draw=Cycle,line width=1.05pt,-{Latex[length=2.8mm,width=1.9mm]},shorten <=5.35mm,shorten >=5.35mm] (p5-C) arc[start angle=-36,end angle=-72,x radius=2.85cm,y radius=3.6000cm];
					\node[anchor=west,font=\small,text=Cycle,fill=white,inner sep=1.2pt] at (6.5368,2.5674) {$CCD$};
					\draw[draw=Cycle,line width=1.05pt,-{Latex[length=2.8mm,width=1.9mm]},shorten <=5.35mm,shorten >=5.35mm] (p5-D) arc[start angle=-72,end angle=-108,x radius=2.85cm,y radius=3.6000cm];
					\node[anchor=north,font=\small,text=Cycle,fill=white,inner sep=1.2pt] at (4.6500,1.8323) {$DDW_1$};
					\draw[draw=Cycle,line width=1.05pt,-{Latex[length=2.8mm,width=1.9mm]},shorten <=5.35mm,shorten >=5.35mm] (p5-Wlast) arc[start angle=-180,end angle=-216,x radius=2.85cm,y radius=3.6000cm];
					\node[anchor=east,font=\small,text=Cycle,fill=white,inner sep=1.2pt] at (1.5971,6.8710) {$W_sW_sX$};
					\draw[draw=Cycle,line width=1.05pt,-{Latex[length=2.8mm,width=1.9mm]},shorten <=5.35mm] (p5-Wone) arc[start angle=-108,end angle=-135,x radius=2.85cm,y radius=3.6000cm];
					\fill[Cycle] (2.4351,3.4159) circle[radius=.035cm];
					\fill[Cycle] (2.2598,3.7208) circle[radius=.035cm];
					\fill[Cycle] (2.1106,4.0471) circle[radius=.035cm];
					\draw[draw=Cycle,line width=1.05pt,-{Latex[length=2.8mm,width=1.9mm]},shorten >=5.35mm] (1.9893,4.3914) arc[start angle=-159,end angle=-180,x radius=2.85cm,y radius=3.6000cm];
					\node[anchor=east,align=right,font=\small,text=Cycle,fill=white,inner sep=1pt] at (1.8184,3.4841) {$W_iW_iW_{i+1}$\\[3pt]$1\leq i<s$};
					\draw[draw=Cycle,line width=1.05pt,-{Latex[length=2.8mm,width=1.9mm]},shorten <=5.35mm,shorten >=5.35mm,draw=Reservoir] (p5-R)--(p5-B);
					\node[font=\small,text=Reservoir] at (8.2750,6.3415) {$RRB$};
					\node[labelled,minimum size=10.0mm] at (p5-X) {$X$};
					\node[labelled,minimum size=10.0mm] at (p5-Y) {$Y$};
					\node[labelled,minimum size=10.0mm] at (p5-Z) {$Z$};
					\node[testpart,minimum size=10.0mm] at (p5-A) {$A$};
					\node[comparepart,minimum size=10.0mm] at (p5-B) {$B$};
					\node[comparecpart,minimum size=10.0mm] at (p5-C) {$C$};
					\node[testpart,minimum size=10.0mm] at (p5-D) {$D$};
					\node[fillerpart,minimum size=10.0mm] at (p5-Wone) {$W_1$};
					\node[fillerpart,minimum size=10.0mm] at (p5-Wlast) {$W_s$};
					\node[reservoirpart,minimum size=10.0mm] at (p5-R) {$R$};
					% Inner boundary has exactly the test-pair colour, width and dash pattern.
					% Its outer stroke edge meets the outer boundary's inner stroke edge:
					% (5mm-1.2pt)+1.25pt/2 = 5mm-1.15pt/2. No white gap.
					\draw[draw=Test,line width=1.25pt,dash pattern=on 4pt off 2.4pt]
					(p5-R) circle[radius={5mm-1.2pt}];
					%\node[font=\normalsize] at (4.6500,0.3161) {$U\longrightarrow V$\quad means\quad $UUV$};
					%\node[align=center,text width=8.6cm,font=\small,text=Muted] at (4.6500,-0.2031) {Two distinct vertices in $U$, one in $V$.};
					\begin{scope}[shift={(9.50,6.6325)}]
						% Item 3, valid for either Construction 2.3 or Construction 6.1.
						% The square stays fixed; only the positions of surrounding text rows change.
						% The plotted square is the parameter domain [0,1)^2; upper/right sides
						% and the diagonal are excluded. The two axes represent the linear order
						% on [0,1), not the circular order used by the arc-test rule.
						\begin{scope}
							% \node[font=\sffamily\bfseries\normalsize,text=Ink] at (3.65,5.08)
							%    {Item 3 in the label plane};
							\node[font=\small,text=Muted] at (3.6500,4.3305)
							{$u\in U,\quad v\in V,\quad b\in B,\quad c\in C$};
							
							% A square of side 4.5 cm, with the origin at (1.2,-.8).
							\path[fill=Compare!10] (1.2,-.8)--(5.7,-.8)--(5.7,3.7)--cycle;
							\path[fill=CompareAlt!10] (1.2,-.8)--(1.2,3.7)--(5.7,3.7)--cycle;
							
							% Upper/right boundaries stand for the excluded value 1.
							\draw[draw=Muted!60,densely dashed,line width=.65pt]
							(1.2,3.7)--(5.7,3.7)--(5.7,-.8);
							\draw[draw=Ink,line width=.8pt,-{Latex[length=2.1mm]}]
							(1.2,-.8)--(6.28,-.8);
							\draw[draw=Ink,line width=.8pt,-{Latex[length=2.1mm]}]
							(1.2,-.8)--(1.2,4.14);
							
							% White clearance along the diagonal emphasizes its exclusion.
							\draw[draw=white,line width=3pt] (1.2,-.8)--(5.7,3.7);
							\draw[draw=Ink!75,densely dashed,line width=.7pt]
							(1.2,-.8)--(5.7,3.7);
							\draw[draw=Ink!70,fill=white,line width=.65pt]
							(1.2,-.8) circle[radius=.045];
							\draw[draw=Ink!70,fill=white,line width=.65pt]
							(5.7,3.7) circle[radius=.045];
							
							\draw[draw=Ink,line width=.65pt] (5.7,-.86)--(5.7,-.74);
							\draw[draw=Ink,line width=.65pt] (1.14,3.7)--(1.26,3.7);
							\node[font=\small,anchor=north east,text=Ink] at (1.12,-.9) {$0$};
							\node[font=\small,anchor=north,text=Ink] at (5.7,-.91) {$1$};
							\node[font=\small,anchor=east,text=Ink] at (1.04,3.7) {$1$};
							\node[font=\normalsize,text=Ink] at (3.5500,-1.3525) {$\lambda(u)$};
							\node[font=\normalsize,text=Ink,rotate=90] at (.49,1.45) {$\lambda(v)$};
							
							\node[font=\large,text=CompareAlt] at (2.64,2.47) {$uvc$};
							\node[font=\small,text=CompareAlt] at (2.64,1.96)
							{$\lambda(u)<\lambda(v)$};
							\node[font=\large,text=Compare] at (4.30,.57) {$uvb$};
							\node[font=\small,text=Compare] at (4.30,.06)
							{$\lambda(v)<\lambda(u)$};
							
							%\draw[draw=Border,line width=.5pt] (0.4500,-1.7625)--(6.8500,-1.7625);
							\node[font=\normalsize,text=Ink] at (3.6500,-2.1590)
							{$(U,V)\in\{(X,Y),(Y,Z),(Z,X)\}$};
							\node[font=\normalsize,text=Compare] at (3.6500,-2.7830)
							{$uvb\in E\quad\Longleftrightarrow\quad\lambda(v)<\lambda(u)$};
							\node[font=\normalsize,text=CompareAlt] at (3.6500,-3.4070)
							{$uvc\in E\quad\Longleftrightarrow\quad\lambda(u)<\lambda(v)$};
						\end{scope}
					\end{scope}
					\draw[draw=Border,line width=.5pt] (9.9500,2.5515)--(16.3500,2.5515);
					\draw[draw=Compare,line width=1.1pt] (10.0000,1.8250)--(11.2500,1.8250);
					\node[anchor=west,font=\normalsize] at (11.6500,1.8250) {Comparison pairs $XB,\ YB,\ ZB$};
					\draw[draw=CompareAlt,line width=1.1pt] (10.0000,1.1625)--(11.2500,1.1625);
					\node[anchor=west,font=\normalsize] at (11.6500,1.1625) {Comparison pairs $XC,\ YC,\ ZC$};
					\draw[draw=Test,line width=1.25pt,dash pattern=on 4pt off 2.4pt] (10.0000,0.5000)--(11.2500,0.5000);
					\node[anchor=west,font=\normalsize] at (11.6500,0.5000) {Test pairs $XA,\ YD$};
					%\node[font=\small,text=Muted] at (13.1500,1.3500) {Pair links are unoriented.};
					%\node[align=center,text width=7cm,font=\small,text=Muted] at (13.1500,0.6000) {All other pairs of distinct parts are ordinary.};
				\end{tikzpicture}%
				}
			\caption{An illustration of Construction~\ref{cons:s-theta}. Here, arrows $U\to V$ denote the triples $UUV$. The unoriented pairs mark the six comparison pairs and two test pairs. The square illustrates Rule~3.}\label{fi:Consruction}
		\end{figure}
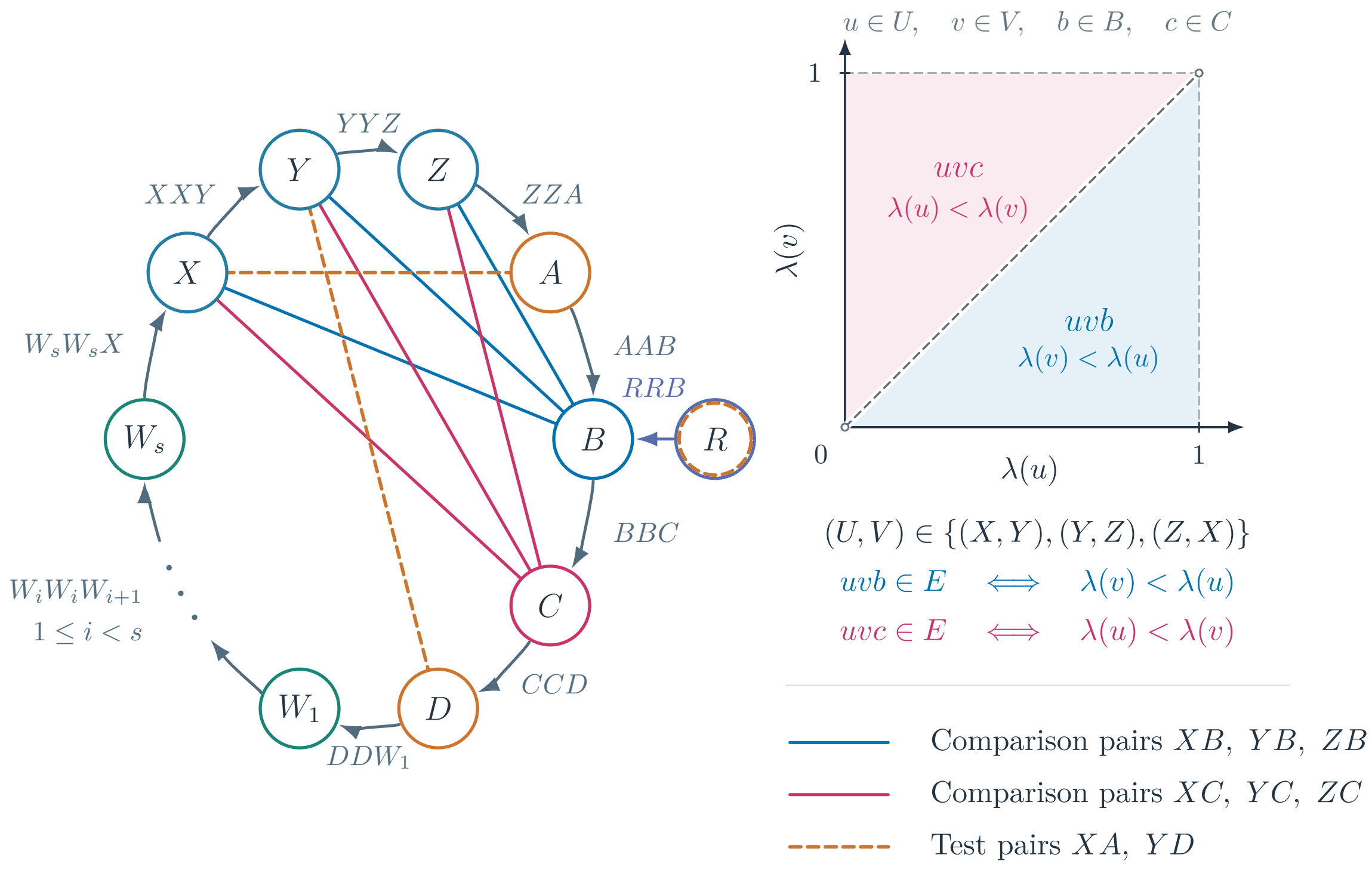	
	\end{center}

	See Figure~\ref{fi:Consruction} illustrating the above definition. Informally speaking, the ideas behind it are as follows. If the minimum codegree is $d\ge2$, a pair within a basic set forces the next set to have at least $d$ vertices. This propagates around the cycle. The comparison pairs constrain the relative orders of the labels in $X,Y,Z$. This allows us to conclude that some test pair receives at most $\theta d$ neighbours from the opposite labelled set, apart from the excess vertices above $d$ in the three labelled sets. The reservoir $R$ must supply the shortfall for the codegree of this test pair, forcing the total number of vertices to be at least $(7+s)d+(1-\theta)d$. Conversely, with evenly spaced labels and basic sets of size $q$, about $(1-\theta)q$ reservoir vertices suffice for codegree $q-1$.

Let us provide all formal details. Note that the defined family $\mathcal A_{s,\theta}$ is closed under taking subgraphs and contains every
	edgeless $3$-graph. Since no part contains an entire edge, $\mathcal A_{s,\theta}$ does not contain
	$K_{2(8+s)+1}^{(3)}$, so $\mathcal F_{s,\theta}$ is nonempty.
	
	We now consider the neighbourhoods of the special pairs in the host defined above. Every cyclically consecutive pair in~\eqref{eq:basic-cycle} is ordinary, as is the pair $\{R,B\}$. Hence, no edge of type $V_iV_iV_{i+1}$ or $RRB$ contains a special pair. The two test pairs consequently have the following neighbourhoods in a host:
	\begin{align}
		N(y,d)
		&=R\cup\{x\in X:\lambda(x)\in I_\theta(\lambda(y))\},
		&&y\in Y,\ d\in D,\label{eq:test-y}\\
		N(x,a)
		&=R\cup\{y\in Y:\lambda(y)\in I_\theta(\lambda(x))\},
		&&x\in X,\ a\in A.\label{eq:test-x}
	\end{align}
	Comparison pairs have neighbours only in the other two labelled classes, as specified in~\eqref{eq:comparison-table}. 
	
	\subsection{The lower bound}
	To prove the lower bound in Theorem~\ref{thm:initial-interval}, for each integer \(n\ge 9+s\) we construct an \(n\)-vertex $3$-graph in \(\mathcal A_{s,\theta}\) with minimum codegree at least \(n/(8+s-\theta)-3\).
	\begin{proposition}\label{prop:base-lower}
		For every integer $s\geq0$, every $\theta\in[0,1]$ and every integer
		$n\geq9+s$,
		\[
		c_{\mathcal A_{s,\theta}}(n)
		\geq\frac{n}{8+s-\theta}-3.
		\]
	\end{proposition}
	
	\begin{proof}
		We first construct $3$-graphs whose orders are not necessarily equal to \(n\), but whose minimum codegrees are easy to calculate. We then delete a few vertices to obtain a graph on exactly \(n\) vertices.
		
		For each integer \(q\geq2\), we construct a graph \(G_q\in\mathcal A_{s,\theta}\) as follows. Take \(q\) vertices in each basic set and \(\lfloor(1-\theta)q\rfloor\) vertices in \(R\). In each of \(X,Y,Z\), choose the labels so that \(\lambda\) takes each value in \(\{0,1/q,\ldots,(q-1)/q\}\) exactly once. Let \(G_q\) contain all triples allowed by the defining edge rules of \(\mathcal A_{s,\theta}\) for these sets and labels.
		
		We now show that \(\delta_2(G_q)=q-1\). A pair within the same basic set has all \(q\) vertices of the next basic set as neighbours. A pair within \(R\) has all \(q\) vertices of \(B\) as neighbours. Every ordinary pair in distinct sets has a whole basic set of neighbours by Rule~2, so its codegree is at least \(q\). Such a set can be chosen from \(X,Y,Z\) if the pair has a labelled vertex, and from \(A,B,C,D\) otherwise. Only comparison and test pairs remain.
		
		For a comparison pair whose labelled vertex has \(\lambda\)-value \(j/q\), there are \(j\) neighbours in one of the other labelled sets and \(q-1-j\) in the other. Its codegree is therefore \(q-1\).
		
		For a test pair, if \(\theta=0\), its neighbourhood is \(R\) and its codegree is \(q\). If \(\theta>0\), its neighbourhood consists of \(R\) and the vertices in the relevant labelled set whose labels lie in the prescribed open arc. This arc contains \(\lceil\theta q\rceil-1\) grid labels, so its codegree is \(\lfloor(1-\theta)q\rfloor+\lceil\theta q\rceil-1=q-1\). Hence \(\delta_2(G_q)=q-1\).
		
		We now construct a graph on exactly \(n\) vertices for each given \(n\geq9+s\). Take \(q=\lceil n/(8+s-\theta)\rceil\geq2\). Since there are \(7+s\) basic sets,
		\[
		n\leq |V(G_q)|=(7+s)q+\lfloor(1-\theta)q\rfloor=\lfloor(8+s-\theta)q\rfloor\leq n+7+s.
		\]
		We can delete \(|V(G_q)|-n\leq7+s\) vertices from distinct basic sets. Every comparison or test pair loses at most two neighbours, while every other pair still has at least \(q-1\) neighbours. The resulting \(n\)-vertex graph belongs to \(\mathcal A_{s,\theta}\) and has minimum codegree at least
		\[
		q-3\geq\frac{n}{8+s-\theta}-3,
		\]
		as required.
	\end{proof}
	
	\subsection{The upper bound}\label{sec:upper-bound}
	We now prove the upper bound in Theorem~\ref{thm:initial-interval}. 
	In our notation, the argument of Mubayi and Zhao~\cite[proof of Theorem~1.6]{MZ07} shows that for every finite set $V_s$ of size $d>0$, every labelling $\lambda:V_s\to[0,1)$ and every $\theta\in[0,1)$, there is a vertex $v\in V_s$ with $|\{w\in V_s: \lambda(w)\in I_\theta(\lambda(v))\}|\le \theta d$. This was very useful in proving existence results for the codegree spectrum of $3$-graphs. We need a generalisation of this result that bounds the number of vertices $w$ in a different part $V_t$ with $\lambda(w)\in I_\theta(\lambda(v))$. This generalisation will be applied to the parts $X,Y,Z$. Using it, we will find a vertex $v$ in $X$ or $Y$ such that the set $I_\theta(\lambda(v))$ contains few labels from the other part. The corresponding test pair (when we combine $v$ with an element of $A$ or $D$) then gives the desired upper bound on the codegree.
	
	\begin{lemma}\label{lem:key}
		Let $V_0,V_1,V_2$ be finite sets with labels in $[0,1)$, where
		$|V_i|=d+e_i$, $d$ is a positive integer and $e_i\geq0$.
		Put $E\coloneqq e_0+e_1+e_2$, and for any $x\in[0,1)$ define
		\[
		L_i(x)\coloneqq|\{v\in V_i:\lambda(v)<x\}|,\qquad
		U_i(x)\coloneqq|\{v\in V_i:\lambda(v)>x\}|.
		\]
		Suppose that, for all distinct $i,j,\ell\in\{0,1,2\}$ and every
		label $x$ occurring in $V_i$,
		\begin{equation}\label{eq:measure-hyp}
			L_j(x)+U_\ell(x)\geq d.
		\end{equation}
		Then, for every $\theta\in[0,1]$ and distinct indices $s,t$,
		at least one of the following holds:
		\begin{enumerate}
			\item There is a vertex $v\in V_s$ such that
			\[
			|\{w\in V_t:\lambda(w)\in I_\theta(\lambda(v))\}|
			\leq\theta d+E.
			\]
			\item There is a vertex $v\in V_t$ such that
			\[
			|\{w\in V_s:\lambda(w)\in I_\theta(\lambda(v))\}|
			\leq\theta d+E.
			\]
		\end{enumerate}
	\end{lemma}
	
	\begin{proof}
		Fix $\theta$. By symmetry, it suffices to consider $V_1$ and $V_2$.
		Choose $i\in\{1,2\}$ and a label $x$ occurring in $V_i$
		that minimize $dx+U_i(x)-e_i$.
		After interchanging $V_1$ and $V_2$ if necessary, assume that $i=1$.
		Then, for every label $y$ occurring in $V_2$,
		\begin{equation}\label{eq:extremal-label}
			dx+U_1(x)-e_1\leq dy+U_2(y)-e_2.
		\end{equation}
		
		We claim that any vertex of $V_1$ with label $x$ satisfies
		the first conclusion of the lemma. Suppose for a contradiction that
		\[
		q\coloneqq|\{w\in V_2:\lambda(w)\in I_\theta(x)\}|
		>\theta d+E.
		\]
		Since $q>0$, we may choose the last label $y$ of $V_2$
		in $I_\theta(x)$, in the forward order from $x$.
		If $y>x$, then $q=U_2(x)-U_2(y)$.
		If $y<x$, then the arc passes through $0$, and
		$q=U_2(x)+d+e_2-U_2(y)$.
		Thus, by~\eqref{eq:extremal-label} and $e_2\geq0$,
		\[
		q\leq d\bigl((y-x)\bmod1\bigr)
		+e_1+U_2(x)-U_1(x)
		<\theta d+e_1+U_2(x)-U_1(x).
		\]
		Consequently,
		\begin{equation}\label{eq:tail-gap}
			U_2(x)>U_1(x)+e_2+e_0.
		\end{equation}
		By~\eqref{eq:measure-hyp} at $x$ and
		$L_2(x)+U_2(x)\leq d+e_2$, we obtain
		\[
		U_0(x)\geq d-L_2(x)
		\geq U_2(x)-e_2
		>U_1(x)+e_0.
		\]
		Together with~\eqref{eq:tail-gap}, this gives
		\[
		L_j(x)\leq d+e_j-U_j(x)<d-U_1(x)
		\qquad(j\in\{0,2\}).
		\]
		
		Let $z$ be the smallest label occurring in $V_0\cup V_2$
		with $z\geq x$. Such a label exists since $U_2(x)>0$.
		By the choice of $z$, neither $V_0$ nor $V_2$ has a label
		in $[x,z)$. Hence
		\[
		L_j(z)=L_j(x)\quad(j\in\{0,2\}),
		\qquad U_1(z)\leq U_1(x).
		\]
		If $z$ occurs in $V_2$, then
		\[
		L_0(z)+U_1(z)
		\leq L_0(x)+U_1(x)<d,
		\]
		contrary to~\eqref{eq:measure-hyp} at $z$.
		The case where $z$ occurs in $V_0$ is identical.
		This proves the claim and hence the lemma.
	\end{proof}
	
	\begin{proposition}\label{prop:base-upper}
		For every integer $s\geq0$, every $\theta\in[0,1]$ and every integer
		$n\geq9+s$,
		\[
		c_{\mathcal A_{s,\theta}}(n)\leq\frac n{8+s-\theta}.
		\]
	\end{proposition}
	
	\begin{proof}
		Take an $n$-vertex $3$-graph $G\in\mathcal A_{s,\theta}$ and fix one of
		its representations. Put $d\coloneqq\delta_2(G)$ and $m\coloneqq7+s$.
		If $d\leq1$, the bound follows from $n\geq m+2$, so assume $d\geq2$.
		
		Since $n\geq m+2$, some vertex class contains two vertices.
		A pair in a basic class has neighbours only in the next basic class.
		The latter therefore has at least $d$ vertices and itself contains a
		pair. Following the cycle~\eqref{eq:basic-cycle} shows that every basic
		class has size at least $d$. If the initial pair lies in $R$, its
		neighbourhood is contained in $B$, and the same argument starts from $B$.
		
		Write
		\[
		|U|=d+\xi_U\quad(U\text{ basic}),\qquad
		S\coloneqq\sum_{U\text{ basic}}\xi_U,
		\qquad E\coloneqq\xi_X+\xi_Y+\xi_Z,
		\qquad r\coloneqq|R|.
		\]
		All excesses are nonnegative, $E\leq S$, and
		\begin{equation}\label{eq:total-mass}
			n=md+S+r.
		\end{equation}
		In particular, each of $X,Y,Z,A,B,C,D$ has size at least $d$.
		
		For $T\in\{X,Y,Z\}$, define $L_T$ and $U_T$ as in
		Lemma~\ref{lem:key}. Let $(T,V,W)$ be any permutation of $(X,Y,Z)$. Consider a vertex of label
		$x$ in $T$. Its comparison pairs with vertices of $B$ and $C$,
		together with~\eqref{eq:comparison-table}, give
		\[
		L_V(x)+U_W(x)\geq d.
		\]
		Thus $X,Y,Z$ satisfy the hypotheses of Lemma~\ref{lem:key},
		with excesses $\xi_X,\xi_Y,\xi_Z$.
		
		Applying Lemma~\ref{lem:key} with $X$ and $Y$ as the two
		distinguished classes, we obtain either a vertex $y\in Y$ such that
		\[
		|\{x\in X:\lambda(x)\in I_\theta(\lambda(y))\}|
		\leq\theta d+E,
		\]
		or a vertex $x\in X$ such that
		\[
		|\{y\in Y:\lambda(y)\in I_\theta(\lambda(x))\}|
		\leq\theta d+E.
		\]
		In the first case, pair $y$ with any vertex of $D$.
		In the second case, pair $x$ with any vertex of $A$.
		The test neighbourhoods~\eqref{eq:test-y} and~\eqref{eq:test-x}
		give, in either case,
		\[
		d\leq r+\theta d+E\leq r+\theta d+S.
		\]
		Combining this with~\eqref{eq:total-mass}, we obtain
		\[
		n=md+S+r\geq(m+1-\theta)d.
		\]
		The result follows from $m+1-\theta=8+s-\theta$.
	\end{proof}
	
	\begin{proof}[Proof of Theorem~\ref{thm:initial-interval}]
		Apply Propositions~\ref{prop:base-lower} and~\ref{prop:base-upper} and
		use~\eqref{eq:class}. Dividing their bounds by $n$ and letting
		$n\to\infty$ gives
		$\gamma(\mathcal F_{s,\theta})=1/(8+s-\theta)$.
	\end{proof}
	
	\section{Extending the interval of densities}\label{sec:spectrum}
	In this section, we extend the interval in Corollary~\ref{cor:base} to $[0,1)$. We use a cyclic substitution in which each part contains a graph from a prescribed subgraph-closed class. The following construction is rather standard, see e.g.\ \cite[Section~3.3]{MZ07}. 
	
	\begin{construction}[$\mathcal D_a(\mathcal C)$]
		Let $a\geq2$ be an integer, and let $\mathcal C$ be a subgraph-closed
		class containing all edgeless $3$-graphs. Define $\mathcal D_a(\mathcal C)$
		to be the class of $3$-graphs $G$ admitting a partition
		\[
		V(G)=V_0\cup\cdots\cup V_{a-1}
		\]
		with the following properties, where empty parts are allowed and indices are taken modulo $a$.
		\begin{enumerate}
			\item The induced graph $G[V_i]$ belongs to $\mathcal C$ for every $i$.
			\item An edge meeting exactly two parts has two vertices in $V_i$ and
			one in $V_{i+1}$ for some $i$.
			\item Edges meeting three distinct parts are unrestricted.
		\end{enumerate}
		Clearly $\mathcal D_a(\mathcal C)$ is also subgraph-closed.
	\end{construction}
	\begin{proposition}\label{prop:substitution}
		Suppose that $\mathcal C$ is subgraph-closed with $\operatorname{dens}(\mathcal C) =\beta$. Then, for every integer
		$a\geq2$,
		\begin{equation}\label{eq:successor-density}
			\operatorname{dens}(\mathcal D_a(\mathcal C))
			=\max\left\{\beta,\frac{1+\beta}{a}\right\}.
		\end{equation}
		In particular, the density is $(1+\beta)/a$ when
		$\beta\leq1/(a-1)$. If $\mathcal C$ omits $K_t^{(3)}$, then
		$\mathcal D_a(\mathcal C)$ omits $K_{a(t-1)+1}^{(3)}$.
	\end{proposition}
	
	\begin{proof}
		Let $G\in\mathcal D_a(\mathcal C)$ have $n>a$ vertices, with a
		partition $V_0,\ldots,V_{a-1}$ as in the definition. Write
		$n_i:=|V_i|$ and $d:=\delta_2(G)$.
		For each $i$ with $n_i\geq2$, choosing a pair of minimum codegree
		in $G[V_i]$ gives
		\[
		d\leq n_{i+1}+c_{\mathcal C}(n_i)
		\leq n_{i+1}+\beta n_i+o(n).
		\]
		Here the error is uniform for $2\leq n_i\leq n$, since $c_{\mathcal C}(t)=\beta t+o(t)$ as $t\to\infty$. If every $n_i\geq2$, summing over $i$ yields
		\[
		d\leq\frac{1+\beta}{a}\,n+o(n).
		\]
		Otherwise, since $n>a$, some index $i$ satisfies
		$n_i\geq2$ and $n_{i+1}\leq1$, giving
		$d\leq\beta n+o(n)$. This proves the upper bound.
		
		For the lower bound, putting all vertices in one part gives
		$\mathcal C\subseteq\mathcal D_a(\mathcal C)$.
		Alternatively, partition $n$ vertices into $a$ parts whose sizes
		differ by at most one, place an extremal member of $\mathcal C$
		in each part, and include all permitted edges between parts.
		Pairs within a part have at least
		$(1+\beta)n/a-o(n)$ neighbours, while pairs in different parts
		have at least $2n/a-O(1)$ neighbours.
		Since $\beta\leq1$, this gives the lower bound $(1+\beta)/a$.
		
		Finally, any copy of $K_{a(t-1)+1}^{(3)}$ has at least $t$ vertices
		in one part, giving a copy of $K_t^{(3)}$ in a member of
		$\mathcal C$, a contradiction.
	\end{proof}
	
	Let $S$ be the set of codegree densities of subgraph-closed classes that omit a complete $3$-graph.
	Corollary~\ref{cor:base} gives $[0,1/7)\subseteq S$. By
	Proposition~\ref{prop:substitution},
	\begin{equation}\label{eq:closure}
		\beta\in S,\qquad 0\leq\beta\leq\frac1{a-1}
		\quad\Longrightarrow\quad
		f_a(\beta):=\frac{1+\beta}{a}\in S
		\qquad(a\geq2).
	\end{equation}
	
	\begin{lemma}\label{lem:interval-extension}
		If $a\geq2$ is an integer and $[0,1/a)\subseteq S$, then
		$[0,1/(a-1))\subseteq S$.
	\end{lemma}
	
	\begin{proof}
		For each integer $j\geq0$, put
		\begin{equation}\label{eq:interval-endpoints}
			r_j:=\frac{1-a^{-j}}{a-1}.
		\end{equation}
		Then $r_0=0$, $r_1=1/a$, and $f_a(r_j)=r_{j+1}$.
		By repeated application of~\eqref{eq:closure},
		\[
		f_a^j([0,1/a))=[r_j,r_{j+1})\subseteq S.
		\]
		Since $r_j\uparrow1/(a-1)$, taking the union over $j\geq0$
		proves the claim.
	\end{proof}
	
	We record the following consequence for use in
	Section~\ref{sec:finite-gaps}.
	
	\begin{proposition}\label{prop:realization}
		For every $\alpha\in[0,1)$, there are $\beta\in[0,1/7)$ and a finite
		sequence $a_1,\ldots,a_t\in\{2,\ldots,7\}$ such that, for every
		subgraph-closed class $\mathcal C$ of density $\beta$ containing all
		edgeless $3$-graphs and omitting a complete $3$-graph, the class
		\[
		\mathcal T(\mathcal C):=
		\mathcal D_{a_t}\bigl(\cdots\mathcal D_{a_1}(\mathcal C)\cdots\bigr)
		\]
		has density $\alpha$ and omits a complete $3$-graph.
		The empty sequence is allowed, with
		$\mathcal T(\mathcal C)=\mathcal C$.
	\end{proposition}
	
	\begin{proof}
		Apply Lemma~\ref{lem:interval-extension} successively for
		$a=7,6,\ldots,2$. Its proof expresses every $\alpha\in[0,1)$ as
		\[
		\alpha=f_{a_t}\circ\cdots\circ f_{a_1}(\beta)
		\]
		for some $\beta\in[0,1/7)$, with each step in the range
		required by~\eqref{eq:closure}.
		By Proposition~\ref{prop:substitution}, the corresponding
		substitutions give density $\alpha$ for every input class of
		density $\beta$ and preserve the omission of a complete $3$-graph.
	\end{proof}
	
	\section{Strict inequalities for finite subfamilies}\label{sec:finite-gaps}
	It remains to prove the finite-subfamily part in Theorem~\ref{thm:spectrum}. 
	We will first prove Lemma~\ref{lem:decreasing-approximation} below, which says that it suffices to express our target class $\mathcal C$ as the intersection of a decreasing sequence of larger classes, each having density strictly greater than that of $\mathcal C$; this idea appears in~\cite{LinZhou2026}.
	We then construct such classes for densities in $[0,1/7)$ which, by cyclic substitution, suffices to cover all densities in $[0,1)$.
	
	\subsection{Decreasing sequences}
	
	\begin{lemma}\label{lem:decreasing-approximation}
		Let $\mathcal C$ be a subgraph-closed class containing every edgeless
		$3$-graph, with $\operatorname{dens}(\mathcal C)=\beta$. Suppose that
		there are subgraph-closed classes
		$\mathcal C_1\supseteq\mathcal C_2\supseteq\cdots$, each containing
		every edgeless $3$-graph, such that
		\begin{equation}\label{eq:decreasing-approximation}
			\mathcal C=\bigcap_{j\geq1}\mathcal C_j,
			\quad \text{and}\quad 
			\beta_j\coloneqq\operatorname{dens}(\mathcal C_j)>\beta
			\quad\text{for every }j\geq1.
		\end{equation}
		Then
		\begin{equation}\label{eq:strict-finite-gap}
			\gamma(\mathcal F')>\beta
			\qquad\text{for every finite }
			\mathcal F'\subseteq\mathcal F(\mathcal C).
		\end{equation}
	\end{lemma}
	
	\begin{proof}
		Let $\mathcal F'$ be a finite subfamily of $\mathcal F(\mathcal C)$.
		If $\mathcal F'=\varnothing$, then
		$\gamma(\mathcal F')=1>\beta$.
		We may therefore assume that $\mathcal F'$ is nonempty.
		
		For every $F\in\mathcal F'$, we have $F\notin\mathcal C$.
		Since $\mathcal C=\bigcap_{j\geq1}\mathcal C_j$, there is an index
		$j_F$ such that $F\notin\mathcal C_{j_F}$. Set
		$j\coloneqq\max_{F\in\mathcal F'}j_F$. By the decreasing property of the sequence,
		we conclude that no member of $\mathcal F'$ belongs to $\mathcal C_j$. Since $\mathcal{C}_j$ is subgraph-closed, every graph in $\mathcal C_j$ is $\mathcal F'$-free.
		Consequently,
		\[
		\gamma(\mathcal F')
		\geq\operatorname{dens}(\mathcal C_j)
		=\beta_j>\beta.
		\]
	\end{proof}
	
	The next corollary shows that cyclic substitution preserves the properties of Lemma~\ref{lem:decreasing-approximation}.
	
	\begin{corollary}\label{cor}
		Let $\mathcal C$ and $(\mathcal C_j)_{j\geq1}$ satisfy the assumptions of
		Lemma~\ref{lem:decreasing-approximation}.
		Then, for every integer $a\geq2$, the classes
		$\mathcal D_a(\mathcal C_j)$ form a decreasing sequence satisfying
		\[
		\mathcal D_a(\mathcal C)=\bigcap_{j\geq1}\mathcal D_a(\mathcal C_j)
		\quad\text{and}\quad
		\operatorname{dens}(\mathcal D_a(\mathcal C_j))
		>\operatorname{dens}(\mathcal D_a(\mathcal C))
		\quad\text{for every }j\geq1.
		\]
		Consequently,
		\[
		\gamma(\mathcal F')>\operatorname{dens}(\mathcal D_a(\mathcal C))
		\qquad\text{for every finite }
		\mathcal F'\subseteq\mathcal F(\mathcal D_a(\mathcal C)).
		\]
	\end{corollary}
	
	\begin{proof}
		From the definition of $\mathcal D_a(\mathcal C_j)$, we know that the classes $\mathcal D_a(\mathcal C_j)$ are subgraph-closed and contain every edgeless $3$-graph. Since the sequence $(\mathcal C_j)_{j\geq1}$ is decreasing, it follows that
		\[
		\mathcal D_a(\mathcal C_1)\supseteq \mathcal D_a(\mathcal C_2)\supseteq\cdots,
		\] 
		forming a decreasing sequence as well.
		
		We next show that their intersection is $\mathcal D_a(\mathcal C)$, that is,
		\begin{equation}\label{eq:successor-intersection}
			\mathcal D_a(\mathcal C)
			=\bigcap_{j\geq1}\mathcal D_a(\mathcal C_j).
		\end{equation}
		For the forward inclusion, let $G\in\mathcal D_a(\mathcal C)$
		and take a partition representing it.
		Since $\mathcal C\subseteq\mathcal C_j$, each part induces a
		$3$-graph in $\mathcal C_j$, while the restrictions on edges
		between parts remain unchanged.
		Thus the same partition shows that
		$G\in\mathcal D_a(\mathcal C_j)$ for every $j\geq1$.
		For the reverse inclusion, fix a finite $3$-graph
		$G\in\bigcap_{j\geq1}\mathcal D_a(\mathcal C_j)$.
		For each $j$, let $\mathcal P_j$ be the set of ordered partitions
		of $V(G)$ into $a$ parts representing $G$ as a member of
		$\mathcal D_a(\mathcal C_j)$.
		Since $\mathcal C_{j+1}\subseteq\mathcal C_j$ and the restrictions
		on edges between parts do not depend on $j$, every partition
		in $\mathcal P_{j+1}$ also belongs to $\mathcal P_j$.
		Thus the sets $\mathcal P_j$ are nonempty, finite and decreasing,
		so some partition belongs to all of them.
		For this partition, each part induces a $3$-graph belonging to
		every $\mathcal C_j$, and hence to $\mathcal C$.
		The same partition therefore represents $G$ as a member of
		$\mathcal D_a(\mathcal C)$, proving
		\eqref{eq:successor-intersection}.
		
		To apply Lemma~\ref{lem:decreasing-approximation}, it remains to show
		that for every $j\ge 1$,
		\[
		\operatorname{dens}(\mathcal D_a(\mathcal C_j))
		>\operatorname{dens}(\mathcal D_a(\mathcal C)).
		\]
		Write $\beta=\operatorname{dens}(\mathcal C)$ and
		$\beta_j=\operatorname{dens}(\mathcal C_j)$, and define
		\[
		g_a(x)\coloneqq\max\left\{x,\frac{1+x}{a}\right\}.
		\]
		By Proposition~\ref{prop:substitution},
		\[
		\operatorname{dens}(\mathcal D_a(\mathcal C_j))=g_a(\beta_j)
		\quad\text{and}\quad
		\operatorname{dens}(\mathcal D_a(\mathcal C))=g_a(\beta).
		\]
		Since $g_a$ is strictly increasing and $\beta_j>\beta$,
		the required inequality follows.
	\end{proof}
	
	\subsection{Small densities}
	
	In this subsection, we construct the decreasing sequences for small densities. For small positive real numbers, we select a suitable sequence of arc parameters in Section~\ref{sec:interval}. The case of density zero, however, requires a separate construction.
	
	\begin{proposition}\label{prop:strict-base}
		For every $\beta\in[0,1/7)$, there is a subgraph-closed class of
		codegree density $\beta$ that contains every edgeless $3$-graph, omits
		a complete $3$-graph, and has a decreasing sequence satisfying
		\eqref{eq:decreasing-approximation}.
	\end{proposition}
	
	\begin{proof}
		Suppose first that $0<\beta<1/7$. Choose an integer $s\geq0$
		and $\theta\in[0,1)$ such that $8+s-\theta=1/\beta$, and set
		\begin{equation}\label{eq:strict-base-parameters}
			t_j\coloneqq\theta+\frac{1-\theta}{j}
			\qquad(j\geq1).
		\end{equation}
		
		Note that $1\ge t_j>\theta$ for every $j$, and $(t_j)_{j\geq1}$ decreases to $\theta$.
		Notice that for every fixed $s$, the classes $\mathcal A_{s,t}$ increase with $t$, as
		increasing the arc length only allows more edges in each $3$-graph. Therefore the classes $(\mathcal A_{s,t_j})_{j\ge 1}$ are decreasing.
		Define
		\begin{equation}\label{eq:strict-base-class}
			\mathcal B_{s,\theta}
			\coloneqq\bigcap_{j\geq1}\mathcal A_{s,t_j}.
		\end{equation}
		We will show that $\mathcal B_{s,\theta}$ has the required properties,
		with $(\mathcal A_{s,t_j})_{j\geq1}$ as its decreasing sequence.
		
		By monotonicity, $\mathcal A_{s,\theta}\subseteq\mathcal A_{s,t_j}$
		for every $j\geq1$. Thus, it holds for every $j\geq 1$ that
		\[
		\mathcal A_{s,\theta}
		\subseteq\mathcal B_{s,\theta}
		\subseteq\mathcal A_{s,t_j}.
		\]
		Then Theorem~\ref{thm:initial-interval} gives, for every integer
		$n\geq9+s$,
		\begin{equation}\label{eq:strict-base-bounds}
			\beta n-3
			\leq c_{\mathcal A_{s,\theta}}(n)\leq c_{\mathcal B_{s,\theta}}(n)\leq \inf_{j\geq1}c_{\mathcal A_{s,t_j}}(n)
			\leq\inf_{j\geq1}\frac{n}{8+s-t_j}
			=\beta n.
		\end{equation}
		It follows that $\operatorname{dens}(\mathcal B_{s,\theta})=\beta$.
		On the other hand, by Theorem~\ref{thm:initial-interval}, we have
		\[
		\operatorname{dens}(\mathcal A_{s,t_j})
		=\frac{1}{8+s-t_j}>\beta
		\qquad(j\geq1).
		\]
		The sequence $(\mathcal A_{s,t_j})_{j\ge 1}$ thus satisfies~\eqref{eq:decreasing-approximation}. Its intersection $\mathcal B_{s,\theta}$ is subgraph-closed and contains every edgeless $3$-graph. It also omits a complete $3$-graph, since
		$\mathcal B_{s,\theta}\subseteq\mathcal A_{s,1}$ and
		$\mathcal A_{s,1}$ has this property by Construction~\ref{cons:s-theta}.
		
		It remains to treat $\beta=0$. Let $\mathcal E$ be the class of
		edgeless $3$-graphs, and let
		\begin{equation}\label{eq:strict-zero-classes}
			\mathcal R_a\coloneqq\mathcal D_a(\mathcal E)
			\quad(a\geq3),
			\qquad
			\mathcal Z_j\coloneqq\bigcap_{a=3}^{j+2}\mathcal R_a,
			\qquad
			\mathcal Z\coloneqq\bigcap_{j\geq1}\mathcal Z_j.
		\end{equation}
		Note that these classes are subgraph-closed and contain every edgeless $3$-graph.
		They also omit a complete $3$-graph. Indeed, $\mathcal R_3$ omits $K_7^{(3)}$, as $\mathcal E$ omits $K_3^{(3)}$. Then all the classes
		$\mathcal Z_j$ and $\mathcal Z$ omit $K_7^{(3)}$ as well.
		By Proposition~\ref{prop:substitution}, we have
		$\operatorname{dens}(\mathcal R_a)=1/a$. It follows that $\operatorname{dens}(\mathcal Z)=0$, as
		$\mathcal Z\subseteq\mathcal R_a$ for every $a\geq3$. 
		
		We now show that $\mathcal Z$ has the required properties,
		with $(\mathcal Z_j)_{j\geq1}$ as its decreasing sequence. The sequence $(\mathcal Z_j)_{j\geq1}$ is decreasing by definition, so it remains to obtain a positive lower bound for the density of each $\mathcal Z_j$. Fix $j$ and let
		\[
		T\coloneqq\prod_{a=3}^{j+2}\mathbb Z/a\mathbb Z,
		\qquad N\coloneqq|T|,
		\]
		where the product is the Cartesian product. Thus, $T$ is the set of
		$j$-tuples $(x_3,\ldots,x_{j+2})$ with
		$x_a\in\mathbb Z/a\mathbb Z$ for each $a\in\{3,\ldots,j+2\}$.
		
		Partition an $n$-vertex set as equally as possible into $N$ parts $B_x$
		indexed by $x\in T$. For three distinct vertices in $B_x,B_y,B_z$, respectively (where $x,y,z$ need not be distinct), include
		their triple if, for every $a\in\{3,\ldots,j+2\}$, the values
		$x_a,y_a,z_a$ are either all distinct or form a multiset
		$\{i,i,i+1\}$ for some $i\in\mathbb Z/a\mathbb Z$.  Denote by $G$ the resulting $3$-graph. We claim that $G$ belongs to $\mathcal Z_j$, namely $G\in \mathcal{R}_a$ for every $a\in\{3,\ldots,j+2\}$.
		Indeed, for each fixed $a\in\{3,\ldots,j+2\}$, grouping the parts $B_x$ together according to the coordinate
		$x_a$ will give a representation in $\mathcal R_a$. Precisely, triples within
		one class are not allowed by our rule, while triples meeting three classes are permitted in $\mathcal R_a$. Moreover, every included triple meeting exactly two classes has two vertices
		in the class indexed by $i$ and one in the class indexed by $i+1$,
		for some $i\in\mathbb Z/a\mathbb Z$.
		
		It remains to bound the minimum codegree of $G$ from below. Consider any pair of distinct vertices in $G$, lying in $B_x$
		and $B_y$, respectively (where $x=y$ is allowed).
		For each $a\in\{3,\ldots,j+2\}$, choose
		\[
		z_a\coloneqq x_a+1\quad\text{if }x_a=y_a,
		\qquad
		z_a\notin\{x_a,y_a\}\quad\text{if }x_a\neq y_a.
		\]
		The second choice is possible as $a\geq3$. The resulting tuple $z$ then differs
		from both $x$ and $y$, and every vertex of $B_z$ forms an edge with the pair. The minimum codegree of $G$ is therefore at least
		$\lfloor n/N\rfloor$, and consequently
		\[
		\operatorname{dens}(\mathcal Z_j)\geq\frac1N>0.
		\]
		Thus the decreasing sequence $(\mathcal Z_j)_{j\geq1}$ satisfies \eqref{eq:decreasing-approximation} with
		$\beta=0$.
	\end{proof}
	
	\subsection{Proof of the main theorem}
	
	\begin{proof}[Proof of Theorem~\ref{thm:spectrum}]
		Fix $\alpha\in[0,1)$. By Proposition~\ref{prop:realization}, there exist $\beta\in[0,1/7)$ and a finite sequence of parameters $a_1,\ldots,a_t$ satisfying the conclusion of Proposition~\ref{prop:realization}. 
		%We list the parameters in the order of application and allow the sequence to be empty.
		
		By Proposition~\ref{prop:strict-base}, choose a subgraph-closed class $\mathcal C_0$
		of density $\beta$ that contains every edgeless $3$-graph, omits a
		complete $3$-graph, and has a decreasing sequence satisfying
		\eqref{eq:decreasing-approximation}. Define
		\[
		\mathcal C_i\coloneqq\mathcal D_{a_i}(\mathcal C_{i-1})
		\qquad(1\leq i\leq t).
		\]
		By Propositions~\ref{prop:realization} and~\ref{prop:substitution}, we have $\operatorname{dens}(\mathcal C_t)=\alpha$, and every $\mathcal C_i$ omits a complete $3$-graph. By repeated applications of
		Corollary~\ref{cor}, the class $\mathcal C_t$
		also has a decreasing sequence satisfying
		\eqref{eq:decreasing-approximation}.
		
		Set $\mathcal F_\alpha\coloneqq\mathcal F(\mathcal C_t)$. Since $\mathcal C_t$ omits a complete $3$-graph and contains
		every edgeless $3$-graph, the family $\mathcal F_\alpha$ is nonempty
		and every member contains an edge. Therefore,
		\[
		\gamma(\mathcal F_\alpha)
		=\operatorname{dens}(\mathcal C_t)=\alpha.
		\]
		Moreover, Lemma~\ref{lem:decreasing-approximation} gives
		\[
		\gamma(\mathcal F')>\alpha
		\qquad\text{for every finite }
		\mathcal F'\subseteq\mathcal F_\alpha.
		\]
		In particular, $\mathcal F_\alpha$ is infinite. This completes the proof.
	\end{proof}
	
	\section{Passing to higher uniformities}\label{sec:lift-uniformity}
	For a $(k+1)$-graph $G$ and $v\in V(G)$, its \emph{link} $L_G(v)$ is the
	$k$-graph with vertex set $V(G)\setminus\{v\}$ and edge set
	\[
	E(L_G(v)):=\{e\subseteq V(G)\setminus\{v\}:e\cup\{v\}\in E(G)\}.
	\]
	Given a family $\mathcal F$ of $k$-graphs, none edgeless, define
	\begin{equation}
		\label{eq:lift}
		\mathcal F^+
		:=
		\left\{G:
		\begin{array}{l}
			G\text{ is a $(k+1)$-graph with }|V(G)|\ge1,\text{ and}\\
			L_G(v)\text{ contains a member of }\mathcal F
			\text{ for every }v\in V(G)
		\end{array}\right\}.
	\end{equation}
	Every member of $\mathcal F^+$ contains an edge.
	For $\mathcal F=\varnothing$, we have $\mathcal F^+=\varnothing$. To prove Theorem~\ref{thm:lift-uniformity}, it suffices to establish
	the following result.
	
	\begin{theorem}
		\label{thm:lift}
		For every integer $k\ge2$ and every family $\mathcal F$ of $k$-graphs,
		\begin{equation}
			\label{eq:equality}
			\gamma(\mathcal F^+)=\gamma(\mathcal F).
		\end{equation}
		Moreover, if $\alpha:=\gamma(\mathcal F)$ and
		$\gamma(\mathcal F_0)>\alpha$ for every finite
		$\mathcal F_0\subseteq\mathcal F$, then
		$\gamma(\mathcal B)>\alpha$ for every finite
		$\mathcal B\subseteq\mathcal F^+$.
	\end{theorem}
	
	Given an ordering $v_1,\ldots,v_n$ of the vertices of a $(k+1)$-graph $G$, for each integer $i\in [n]$ we define its \emph{forward link} $L_i^+(G)$, where
	\[
	\begin{aligned}
		V(L_i^+(G))&:=\{v_{i+1},\ldots,v_n\},\\
		E(L_i^+(G))&:=
		\left\{e\in\binom{\{v_{i+1},\ldots,v_n\}}{k}:
		e\cup\{v_i\}\in E(G)\right\}.
	\end{aligned}
	\]
	
	\begin{lemma}
		\label{lem:order}
		A $(k+1)$-graph $G$ is $\mathcal F^+$-free if and only if its vertices
		admit an ordering for which every forward link is $\mathcal F$-free.
	\end{lemma}
	
	\begin{proof}
		Suppose first that $G$ is $\mathcal F^+$-free.
		We claim that every nonempty induced subgraph $G'$ of $G$ has a vertex $v$
		whose link $L_{G'}(v)$ is $\mathcal F$-free.
		Indeed, otherwise $G'$ belongs to $\mathcal F^+$ by~\eqref{eq:lift}.
		Repeatedly delete such a vertex, and list the vertices in the order
		of deletion. At each step the current link is exactly the forward
		link in the resulting ordering.
		
		Conversely, suppose that such an ordering exists and that $G$
		contains a copy of some $B\in\mathcal F^+$.
		Let $v_i$ be the earliest vertex of that copy.
		The link of $v_i$ in $B$ is a subgraph of
		$L_i^+(G)$. By~\eqref{eq:lift}, this link contains a member of
		$\mathcal F$, contrary to the $\mathcal F$-freeness of $L_i^+(G)$.
	\end{proof}
	
	\begin{proof}[Proof of Theorem~\ref{thm:lift}]
		Let $\mathcal F$ be a family of $k$-graphs with
		$\gamma(\mathcal F)=\alpha$.
		For the upper bound on $\gamma(\mathcal F^+)$, let $G$ be an $\mathcal F^+$-free
		$(k+1)$-graph on $n\geq k+1$ vertices.
		By Lemma~\ref{lem:order}, some vertex $v$ has an
		$\mathcal F$-free link $L_G(v)$. Hence
		\[
		\delta_k(G)
		\leq \delta_{k-1}(L_G(v))
		\leq \operatorname{coex}(n-1,\mathcal F).
		\]
		Thus
		\begin{equation}\label{eq:upper}
			\operatorname{coex}(n,\mathcal F^+)
			\leq \operatorname{coex}(n-1,\mathcal F),
		\end{equation}
		which gives $\gamma(\mathcal F^+)\leq\alpha$.
		
		For the lower bound, we may assume that $\alpha>0$.
		Fix $0<\varepsilon<\alpha$ and put $\beta:=\alpha-\varepsilon$.
		There is an integer $M\geq k$ such that, for every $m\geq M$,
		we can choose an $\mathcal F$-free $k$-graph $H_m$ on $m$
		vertices with
		\begin{equation}\label{eq:input}
			\delta_{k-1}(H_m)\geq\beta m.
		\end{equation}
		Double counting gives
		\begin{equation}\label{eq:edge-density}
			\frac{e(H_m)}{\binom{m}{k}}
			\geq \frac{\delta_{k-1}(H_m)}{m-k+1}
			\geq\beta.
		\end{equation}
		
		Let $n>M$ be a sufficiently large integer. For each $1\leq i\leq n-M$, independently choose
		a uniformly random bijection from $V(H_{n-i})$ to
		$\{i+1,\ldots,n\}$, and let $J_i$ be the resulting labelled
		copy of $H_{n-i}$.
		For $n-M<i\leq n$, let $J_i$ be the edgeless $k$-graph
		on $\{i+1,\ldots,n\}$.
		Define a random $(k+1)$-graph $G$ on $[n]$ by
		\begin{equation}\label{eq:construction}
			E(G):=\bigcup_{i=1}^n
			\{\{i\}\cup e:e\in E(J_i)\}.
		\end{equation}
		By construction, the forward link at each vertex $i$ is
		exactly $J_i$, which is $\mathcal F$-free.
		Thus $G$ is $\mathcal F^+$-free by Lemma~\ref{lem:order}.
		
		Fix $S\in\binom{[n]}k$ and put $i:=\min S$.
		Let $X_S$ and $Y_S$ count the neighbours of $S$ smaller
		and larger than $i$, respectively.
		Then $\deg_G(S)=X_S+Y_S$, and
		\[
		Y_S=\deg_{J_i}(S\setminus\{i\}).
		\]
		By~\eqref{eq:input},
		\begin{equation}\label{eq:later}
			Y_S\geq
			\begin{cases}
				\beta(n-i),& i\leq n-M,\\
				0,& i>n-M.
			\end{cases}
		\end{equation}
		
		Set $q_i:=\min\{i-1,n-M\}$.
		By~\eqref{eq:construction},
		\[
		X_S=\sum_{j=1}^{q_i}
		{\bf 1}_{\{S\in E(J_j)\}}.
		\]
		These indicators are independent, and uniform relabelling
		together with~\eqref{eq:edge-density} gives
		\begin{equation}\label{eq:probability}
			\mathbb P(S\in E(J_j))
			=\frac{e(H_{n-j})}{\binom{n-j}{k}}
			\geq\beta.
		\end{equation}
		Thus $\beta q_i\leq\mathbb E X_S\leq n$, and the
		Chernoff bound~\cite[Theorem~2.8]{JLR00} gives
		\[
		\mathbb P(X_S<\beta q_i-\varepsilon n)
		\leq \exp(-\varepsilon^2n/2),
		\]
		where the estimate is trivial if $\beta q_i\leq\varepsilon n$.
		A union bound over all $k$-sets shows that, with probability
		$1-o(1)$,
		\[
		X_S\geq\beta q_i-\varepsilon n
		\qquad\text{for every }S\in\binom{[n]}k.
		\]
		For sufficiently large $n$, choose $G$ with this property. If $i\le n-M$, then $q_i=i-1$ and
		\eqref{eq:later} gives
		\[
		\deg_G(S)\ge\beta(n-i)+\beta(i-1)-\varepsilon n
		=\beta(n-1)-\varepsilon n.
		\]
		If $i>n-M$, then $q_i=n-M$, so
		$\deg_G(S)\ge\beta(n-M)-\varepsilon n$.
		Since $M\ge1$, both cases imply
		\[
		\delta_k(G)\geq\beta(n-M)-\varepsilon n
		=(\alpha-2\varepsilon)n-\beta M.
		\]
		Dividing by $n$ and letting $n\to\infty$ gives
		$\gamma(\mathcal F^+)\geq\alpha-2\varepsilon$.
		Letting $\varepsilon\to0$ and using~\eqref{eq:upper}
		proves~\eqref{eq:equality}.
		
		Suppose now that every finite subfamily of $\mathcal F$
		has density strictly greater than $\alpha$.
		Let $\mathcal B\subseteq\mathcal F^+$ be finite.
		For each $B\in\mathcal B$ and $v\in V(B)$, choose
		$F_{B,v}\in\mathcal F$ that embeds in $L_B(v)$, and put
		\[
		\mathcal F_0:=
		\{F_{B,v}:B\in\mathcal B,\ v\in V(B)\}.
		\]
		Then $\mathcal F_0$ is finite and
		$\mathcal B\subseteq\mathcal F_0^+$. Consequently,
		\[
		\gamma(\mathcal B)
		\geq\gamma(\mathcal F_0^+)
		=\gamma(\mathcal F_0)
		>\alpha,
		\]
		where the equality follows from~\eqref{eq:equality}.
	\end{proof}
	
	\section{Terminal intervals for positive codegree densities}
	\label{sec:positive-terminal}
	
	In this section, we prove Theorem~\ref{thm:positive-introduction}, which states that the positive codegree spectrum of $k$-graphs contains an interval $[c,1)$ for every $k\ge3$. Our goal is to establish the existence of such a constant $c<1$,
	rather than to optimise its value. 
	
	Throughout this section, all
	$k$-graphs are finite. Given a $k$-graph $G$, we call a vertex subset $S$ of $G$ \emph{covered} if $S$ is contained in an edge of $G$. A covered $(k-1)$-set is called a \emph{facet}. Define
	\[
	d(G)\coloneqq \delta_{k-1}^{+}(G)
	=\min\{|N_G(S)|: S\text{ is a facet of }G\},
	\]
	with $d(G)=0$ when $G$ is edgeless. For a subgraph-closed class $\mathcal C$ containing
	all edgeless $k$-graphs, let
	\[
	c^+_{\mathcal C}(n)\coloneqq \max\{d(G):G\in\mathcal C,\ |V(G)|=n\}.
	\]
	As before, $\mathcal C$ is exactly the class of $\mathcal F(\mathcal C)$-free
	graphs. Hence, 
	\begin{equation}\label{eq:positive-class}
		\operatorname{ex}^+(n,\mathcal F(\mathcal C))=c^+_{\mathcal C}(n),
	\end{equation}
	and the limit of $c^+_{\mathcal C}(n)/n$ as $n\to\infty$ exists by~\cite[Proposition~1]{Pikhurko23}.
	
	In this section, we will first build a
	$3$-uniform core family and prove an upper bound under an additional assumption. Then we will use it to obtain
	the desired interval in every uniformity $k\geq3$.
	
	\subsection{A three-uniform core family}
	
	Recall that for
	$t\in[0,1)$ and $\theta\in[0,1]$, the forward arc is defined as
	\[
	I_\theta(t)=\{u\in[0,1):0<(u-t)\bmod1<\theta\}.
	\]
	The following construction is very similar to
	Construction~\ref{cons:s-theta} with $s=1$ and $W_1=W$. The difference is that Construction~\ref{cons:s-theta} allows all triples meeting three distinct ordinary parts while Construction~\ref{cons:pos-core} includes such a  triple only if it intersects $W$ and an unprotected pair of parts. Without this change, Lemma~\ref{lem:pos-core} would fail. For example, a blowup $H$ of $K_4^{(3)}$ on parts $A,B,C,D$ of equal size $n/4$ (with all other parts empty) has  $d(H)=n/2$ and $h(H)=n/4$.

	\begin{construction}[$\mathcal C_\theta$]\label{cons:pos-core}
		Fix $\theta\in[0,1]$ and take pairwise disjoint vertex sets
		\begin{equation}\label{eq:pos-sets}
			X,Y,Z,A,B,C,D,W,R.  
		\end{equation}
		We call the sets other than $R$ \emph{basic sets} and arrange them in the
		cyclic order
		\begin{equation}\label{eq:pos-basic-cycle}
			(V_0,\ldots,V_7)\coloneqq(X,Y,Z,A,B,C,D,W).
		\end{equation}
		Assign each vertex $v\in X\cup Y\cup Z$ a label $\lambda(v)\in [0,1)$.
		Empty sets and repeated labels are allowed.
		
		As before, the six pairs $\{U,B\}$ and $\{U,C\}$ with
		$U\in\{X,Y,Z\}$ are called \emph{comparison pairs}, and
		$\{X,A\}$ and $\{Y,D\}$ are called \emph{test pairs}.
		These eight pairs are \emph{special}; all remaining pairs of distinct
		sets are \emph{ordinary}.
		For the present construction, call a pair \emph{protected} if it is
		special or is one of
		\[
		\{X,Y\},\quad\{Y,Z\},\quad\{Z,X\},\quad\{R,X\},\quad\{R,Y\}.
		\]
		
		For each choice of vertex sets and labels, define a host $3$-graph by the
		following edge rules, with no other edges permitted.
		\begin{enumerate}
			\item Include every triple of type $V_iV_iV_{i+1}$, with indices taken
			modulo eight, and every triple of type $RRB$.
			
			\item For every unprotected pair $\{U,U'\}$ of distinct sets other than
			$W$, include every triple of type $WUU'$. We call such triples
			\emph{filler triples}. 
			
			\item For two vertices in distinct labelled sets and a third vertex in
			$B$ or $C$, include their triple if and only if the corresponding
			inequality in the following table holds:
			\begin{equation}\label{eq:pos-comparison-table}
				\begin{array}{c|cc}
					\text{labelled vertices}
					&\text{third vertex in }B&\text{third vertex in }C\\\hline
					x\in X,\ y\in Y&\lambda(y)<\lambda(x)&\lambda(x)<\lambda(y)\\
					y\in Y,\ z\in Z&\lambda(z)<\lambda(y)&\lambda(y)<\lambda(z)\\
					z\in Z,\ x\in X&\lambda(x)<\lambda(z)&\lambda(z)<\lambda(x)
				\end{array}
			\end{equation}
			
			\item Include every triple of type $XAR$ or $YDR$.
			For $x\in X$ and $y\in Y$, also include $xya$ with $a\in A$ if
			$\lambda(y)\in I_\theta(\lambda(x))$, and include $xyd$ with $d\in D$
			if $\lambda(x)\in I_\theta(\lambda(y))$.
		\end{enumerate}
		Let $\mathcal C_\theta$ consist of all subgraphs of these hosts, allowing
		all set sizes and all labellings.	\end{construction}
	
	%Every adjacent pair in~\eqref{eq:pos-basic-cycle} is ordinary, as is
	%$\{R,B\}$. Hence the repeated-set triples do not contain a comparison
	%or test pair. The filler triples also avoid all protected pairs.
	Note that the class $\mathcal C_\theta$ is closed under
	taking subgraphs and contains every edgeless $3$-graph.
		The two test pairs have the following neighbourhoods in a host:
	\begin{align}
		N(x,a)&=R\cup\{y\in Y:\lambda(y)\in I_\theta(\lambda(x))\},
		&&x\in X,\ a\in A,\label{eq:pos-test-x}\\
		N(y,d)&=R\cup\{x\in X:\lambda(x)\in I_\theta(\lambda(y))\},
		&&y\in Y,\ d\in D.\label{eq:pos-test-y}
	\end{align}
	Comparison pairs have neighbours only in the other two labelled sets,
	as specified in~\eqref{eq:pos-comparison-table}.
	%In subgraphs, these descriptions give upper bounds on the corresponding
	%neighbourhoods. The additional protected pairs will allow us to force
	%all eight basic sets to be large.
	
	For a nonempty $3$-graph $H$ without isolated vertices, let $\partial_2H$ denote its shadow graph and let $h(H)\coloneqq |V(H)|-\delta(\partial_2H)$. Thus $h(H)$ is the
	largest number of non-neighbours of a vertex in the shadow graph, counting the vertex itself. 
	
	In the following lemma, we are able to provide an upper bound on $d(H)$ for $H\in\mathcal C_\theta$, under an additional assumption.
	
	\begin{lemma}\label{lem:pos-core}
		Set $L\coloneqq 9-\theta$. Let $H\in\mathcal C_\theta$ be a nonempty $3$-graph without isolated vertices. If $H$ satisfies $d(H)>h(H)$, then 
		\[
		d(H)\le |V(H)|/L.
		\]
		Moreover, for every sufficiently large integer $n$, the class $\mathcal C_\theta$ contains an $n$-vertex $3$-graph with minimum codegree at least $n/L-2$.
	\end{lemma}
	\begin{proof}
		Suppose $d=d(H)>h(H)=h$. By the definition of $h$, for every vertex $v$, any set of more than $h$ vertices
		contains a neighbour of $v$ in the shadow graph. In particular, any
		part of size greater than $h$ contains a covered pair. We first prove the following claim.
		
		\begin{claim}\label{cl:positive-set}
			All basic sets in~(\ref{eq:pos-basic-cycle}) have size at least $d$.
		\end{claim}
		
		\begin{poc}
			We first prove that there exists a covered pair lying inside one of the nine sets. Suppose otherwise. We will show that $H$ has no edges, contradicting
			the assumption that $H$ is nonempty. Under this assumption, each of the nine sets has size at most $h<d$, and all triples of type $V_iV_iV_{i+1}$ and $RRB$ are absent.
			
			An edge of type $XAR$ would contain a covered pair with one vertex
			in $R$ and the other in $X$. All neighbours of this pair lie in $A$,
			contradicting $|A|<d$. Similarly, an edge of type $YDR$ would give
			a covered pair between $R$ and $Y$ with all its neighbours in $D$.
			Thus $H$ contains no edge of either type.
			
			It follows that pairs of type $XA$ have neighbours only in $Y$,
			and pairs of type $YD$ have neighbours only in $X$.
			Since $|X|,|Y|<d$, this excludes all edges of types $XYA$ and $XYD$.
			Now consider a pair of vertices from two distinct labelled sets.
			By~\eqref{eq:pos-comparison-table}, the pair has neighbours in at
			most one of $B$ and $C$, depending on the order of its labels.
			Since $|B|,|C|<d$, this excludes all triples containing comparison pairs.
			Only filler triples remain. If $H$ contained a filler triple of
			type $WUU'$, the pair of vertices in $U$ and $U'$ would be covered
			and all its neighbours would lie in $W$, contradicting $|W|<d$.
			Thus $H$ has no edges, contrary to our assumption.
			
			Hence one of the nine sets contains a covered pair. If this pair
			lies in a basic set $V_i$, all its neighbours lie in $V_{i+1}$.
			Therefore $|V_{i+1}|\geq d>h$, so $V_{i+1}$ also contains a covered
			pair. Continuing around the cycle gives $|V_i|\geq d$ for every $i$.
			If the initial covered pair lies in $R$, all its neighbours lie in $B$.
			Thus $|B|\geq d>h$, and we can apply the same argument starting
			with a covered pair in $B$.
		\end{poc}
		
		Set
		\[
		E\coloneqq |X|+|Y|+|Z|-3d.
		\]
		By Claim~\ref{cl:positive-set}, we have $E\geq0$.
		For $T\in\{X,Y,Z\}$, define $L_T$ and $U_T$ as in
		Lemma~\ref{lem:key}:
		\[
		L_T(x)\coloneqq|\{v\in T:\lambda(v)<x\}|,\qquad
		U_T(x)\coloneqq|\{v\in T:\lambda(v)>x\}|.
		\]
		Fix a vertex $v$ in a labelled set $T$, and write $t=\lambda(v)$.
		Since $|B|,|C|\geq d>h$, the vertex $v$ has a neighbour in each
		of $B$ and $C$ in the shadow graph. Each of these two covered pairs has codegree at least $d$. Applying
		\eqref{eq:pos-comparison-table} to these two pairs gives
		\[
		L_{T'}(t)+U_{T''}(t)\geq d
		\]
		for both orderings $(T',T'')$ of the other two labelled sets.
		Thus the hypotheses of Lemma~\ref{lem:key} hold, with total excess $E$.

		Applying Lemma~\ref{lem:key} with $X$ and $Y$ as the distinguished
		classes, we obtain either a vertex of $X$ with at most $\theta d+E$
		vertices of $Y$ in its forward $\theta$-arc, or a vertex of $Y$
		with at most $\theta d+E$ vertices of $X$ in its forward
		$\theta$-arc. In the first case, let $x\in X$ be such a vertex.
		Since $|A|\geq d>h$, there is a vertex $a\in A$ adjacent to $x$
		in the shadow graph. By~\eqref{eq:pos-test-x},
		\[
		d\leq |N_H(x,a)|\leq |R|+\theta d+E.
		\]
		In the second case, choosing a neighbour in $D$ and
		applying~\eqref{eq:pos-test-y} gives the same bound.
		
		By Claim~\ref{cl:positive-set}, the other five basic sets each have
		size at least $d$. Consequently,
		\[
		|V(H)|\geq |X|+|Y|+|Z|+5d+|R|
		=8d+E+|R|
		\geq(9-\theta)d,
		\]
		thus proving the first conclusion of the lemma.

		For the lower bound, fix a positive integer $q$. Put $q$ vertices in each basic set and
		$\lfloor(1-\theta)q\rfloor$ vertices in $R$.
		Write $X=\{x_0,\ldots,x_{q-1}\}$, and similarly for $Y,Z$, and set
		\[
		\lambda(x_j)=\frac{3j}{3q},\qquad
		\lambda(y_j)=\frac{3j+1}{3q},\qquad
		\lambda(z_j)=\frac{3j+2}{3q}
		\quad(0\le j<q).
		\]
		Include every permitted triple. We claim that the resulting $3$-graph
		has minimum codegree at least $q-1$.
		
		First, every pair inside a basic set has all vertices of the next set in
		the cycle as neighbours, and every pair inside $R$ has all of $B$
		as neighbours. An unprotected pair between distinct sets other than
		$W$ has all of $W$ as neighbours. For a pair with exactly one vertex
		in $W$, one of the filler types
		\[
		WXD,\quad WYA,\quad WZA,\quad WRA,\quad WBC
		\]
		provides a full basic set of neighbours.
		
		For a pair in two distinct labelled sets, the two labels are distinct,
		so all vertices of either $B$ or $C$ are neighbours, according
		to~\eqref{eq:pos-comparison-table}. For a comparison pair whose
		labelled vertex has index $j$, the neighbourhood contains all vertices
		with smaller labels in one of the other two labelled sets and all
		vertices with larger labels in the other. There are at least $j$ of
		the former and at least $q-1-j$ of the latter, giving at least $q-1$
		neighbours.
		
		Pairs of type $RX$ and $RY$ have all vertices of $A$ and $D$,
		respectively, as neighbours. For a test pair of type $XA$ or $YD$,
		all vertices of $R$ are neighbours. Moreover, the labels in each
		labelled set are equally spaced around the circle, so the forward
		$\theta$-arc contains at least $\lceil\theta q\rceil-1$ labels
		from the other set ($Y$ or $X$, respectively). Thus each test pair
		has at least
		\[
		\lfloor(1-\theta)q\rfloor+\lceil\theta q\rceil-1=q-1
		\]
		neighbours. This accounts for all pairs. The number of vertices is
		$8q+\lfloor(1-\theta)q\rfloor=\lfloor Lq\rfloor$.
		
		For a sufficiently large integer $n$, take $q=\lfloor n/L\rfloor$
		and add $n-\lfloor Lq\rfloor$ new vertices to $R$, again including every
		permitted triple. The codegree of each pair of old vertices does not
		decrease. Every pair containing a new vertex has a full basic set
		of neighbours: $B$ for a pair inside $R$, $A$ for a pair of type
		$RX$ or $RW$, $D$ for a pair of type $RY$, and $W$ in all remaining
		cases. Hence the resulting $n$-vertex $3$-graph has minimum codegree
		at least $q-1\geq n/L-2$.
	\end{proof}
	
	As $\theta$ ranges over $[0,1]$, the parameter $\alpha=1/(9-\theta)$
	ranges over $[1/9,1/8]$. The lemma does not by itself show that
	$\mathcal C_\theta$ has positive codegree density $\alpha$, since
	its upper bound requires $d(H)>h(H)$.
	
	\subsection{Direct amplification in every uniformity}
	
	We first introduce the following lemma, which will be used when constructing the hypergraph family realising the densities in $[c,1)$.

	\begin{lemma}[Cyclic blocks]\label{lem:pos-blocks}
		Let $m\ge 2$ and $p\ge 2^m$ be integers. There are $p$ distinct
		$m$-element blocks $B_1,\ldots,B_p$ on a set of $p$ points such that
		every point lies in exactly $m$ blocks and any two distinct points
		lie together in at most one block.
	\end{lemma}
	
	\begin{proof}
		Take the point set to be $\mathbb Z/p\mathbb Z$. Let
		$B_1\coloneqq\{1,2,4,\ldots,2^{m-1}\}$ and put
		$B_i\coloneqq B_1+(i-1)$ for $2\le i\le p$.
		
		Suppose that two distinct points $x,y$ lie in $B_i$. Then
		$x=i-1+b$ and $y=i-1+c$ for some $b,c\in B_1$, so
		\[
		b-c\equiv x-y\pmod p.
		\]
		Since $|b-c|\le 2^{m-1}-1<p/2$, this congruence determines $b-c$
		as an integer. Interchanging $x$ and $y$ if necessary, assume that
		$b>c$. Since the positive differences between elements of $B_1$ are all distinct, $b-c$ determines both $b$ and $c$. Once $b$
		is known, the shift is given by $i-1=x-b$. It follows that at most
		one index $i$ is possible for the pair $x,y$.
		
		In particular, the $p$ blocks are distinct, since each contains at
		least two points. Finally, a fixed point $x$ lies in $B_i$ precisely
		when $i-1=x-b$ for some $b\in B_1$. These $m$ choices of $b$ give
		$m$ distinct shifts, so $x$ lies in exactly $m$ blocks.
	\end{proof}
	
	\begin{construction}[$\mathcal P_{k,p,\alpha}$]
		\label{cons:pos-amplification}
		Fix $m=72$, integers $p\ge 2^m$ and $k\ge 3$, and
		$\alpha\in[1/9,1/8]$. Set $\theta\coloneqq 9-1/\alpha$.
		Take pairwise disjoint sets
		\[
		V,Q,U_0,U_1,\ldots,U_p
		\]
		with $Q=\varnothing$ when $k=3$. Let $K\coloneqq V\cup Q$ and
		let $H\in\mathcal C_\theta$ be a $3$-graph on $V$.
		
		We call $U_0,\ldots,U_p$ the \emph{outside colour classes} and
		their vertices \emph{outside vertices}. Fix a block system
		$B_1,\ldots,B_p$ on the colours $1,\ldots,p$ as in
		Lemma~\ref{lem:pos-blocks}. For each $i\in[p]$, we call the vertex set
		\[
		U_{B_i}\coloneqq\bigcup_{x\in B_i}U_x
		\]
		the \emph{colour block} corresponding to $B_i$. Since the colour $0$
		belongs to no block, $U_0$ lies in no colour block.
		
		The host is the $k$-graph on $K\cup U_0\cup\cdots\cup U_p$ whose
		edges are precisely the $k$-sets $E$ satisfying all four of the
		following conditions.
		\begin{enumerate}[label=(\arabic*), ref=(\arabic*)]
			\item\label{item:pos-ampl-core}
			$|E\cap V|\le3$, with $E\cap V\in E(H)$ if equality holds.
			\item\label{item:pos-ampl-internal}
			If $E\subseteq K$, then $|E\cap V|\in\{1,3\}$.
			\item\label{item:pos-ampl-colour}
			Each outside colour class contains at most $k-2$ vertices of $E$.
			\item\label{item:pos-ampl-block}
			No colour block contains $E$.
		\end{enumerate}
		Let $\mathcal P_{k,p,\alpha}$ consist of all subgraphs of these hosts,
		allowing all part sizes and all choices of $H\in\mathcal C_\theta$.
		This class is subgraph-closed and contains every edgeless $k$-graph.
	\end{construction}

	Note that an edge containing a pair $S\subseteq V$ cannot lie in $S\cup Q$,
	since it would then be contained in $K$ and meet $V$ in exactly two
	vertices, contrary to Rule~\ref{item:pos-ampl-internal}. Similarly,
	for $v\in V$ and $0\le i\le p$, an edge containing $v$ cannot lie in
	$\{v\}\cup U_i$, since it would contain $k-1$ vertices of colour $i$,
	contrary to Rule~\ref{item:pos-ampl-colour}. These properties allow applications of the following lemma. Recall that
	$d(G)=\delta^+_{k-1}(G)$ denotes the minimum positive codegree of $G$.
	
	\begin{lemma}\label{lem:pos-extension}
		Let $G$ be a $k$-graph with at least one edge. Let $T$ be a covered
		set of size less than $k$, and let $I\subseteq V(G)\setminus T$.
		If no edge containing $T$ lies in $T\cup I$, then there are at least
		$d(G)$ vertices $x\in V(G)\setminus(T\cup I)$ such that
		$T\cup\{x\}$ is covered.
		
		In particular, $d(G)\le |V(G)|-|I|$ for every independent set $I$ in
		$G$, and every covered set of size less than $k$ has at least $d(G)$
		extensions.
	\end{lemma}
	
	\begin{proof}
		Among the edges containing $T$, choose $E$ with $|E\cap I|$ as large
		as possible. By hypothesis, there exists a vertex $z\in E\setminus(T\cup I)$.
		Then $F\coloneqq E\setminus\{z\}$ is a facet containing $T$, and $F$ has no neighbour in $I$, since such a neighbour would give an edge increasing $|E\cap I|$. 
		Furthermore, by the definition of $d(G)$, the facet
		$F$ has at least $d(G)$ neighbours. For each such neighbour $x$,
		the edge $F\cup\{x\}$ covers $T\cup\{x\}$.
		
		For the first consequence, let $I$ be independent and take
		$T=\varnothing$. Since $I$ contains no edge, the lemma gives at least
		$d(G)$ vertices outside $I$, so $d(G)\le |V(G)|-|I|$. For the second,
		take $I=\varnothing$, and the
		lemma gives at least $d(G)$ extensions of $T$.
	\end{proof}

	Now we are ready to prove Theorem~\ref{thm:positive-introduction}.
	
	\begin{proof}[Proof of Theorem~\ref{thm:positive-introduction}]
		Fix $m=72$ and an integer $k\ge 3$. For each integer $p\ge 2^m$ and
		each $\alpha\in[1/9,1/8]$, let $\mathcal P_{k,p,\alpha}$ be the class
		in Construction~\ref{cons:pos-amplification}, with $\theta=9-1/\alpha$,
		and put
		\begin{equation}\label{eq:pos-R}
			\rho=
			\begin{cases}
				1+\dfrac{1}{1-\alpha}+\dfrac{p-1}{m}, & k=3,\\[4pt]
				2+\alpha+\dfrac{p}{m}, & k\ge 4.
			\end{cases}
		\end{equation}
		
		For fixed $p$, the parameter $\rho$ is continuous and strictly
		increasing in $\alpha$. At $\alpha=1/9$, both cases of~\eqref{eq:pos-R}
		give $\rho=(p+152)/72$. As $\alpha$ increases from $1/9$ to $1/8$,
		the increase in $\rho$ is $1/56$ when $k=3$ and $1/72$ when $k\ge 4$.
		Thus, in either case, $\rho$ takes every value in
		\[
		I_p\coloneqq
		\left[\frac{p+152}{72},\frac{p+153}{72}\right].
		\]
		The right endpoint of $I_p$ is the left endpoint of $I_{p+1}$, so
		the union of these intervals over all integers $p\ge 2^{72}$ is
		$[(2^{72}+152)/72,\infty)$. Hence every number in $[c,1)$ can be
		written as $1-1/\rho$ for some $p$ and $\alpha$, where
		\[
		c\coloneqq 1-\frac{72}{2^{72}+152}.
		\]
		
		By~\eqref{eq:positive-class}, the positive codegree Tur\'an density of
		$\mathcal F(\mathcal P_{k,p,\alpha})$ equals the positive codegree
		density of $\mathcal P_{k,p,\alpha}$. To realise every density in
		$[c,1)$, it therefore suffices to prove that
		\[
		\lim_{n\to\infty}\frac{c^+_{\mathcal P_{k,p,\alpha}}(n)}{n}
		=1-\frac{1}{\rho}.
		\]
		
		\noindent\textbf{The upper bound.}  For the upper bound, we show that every $G\in\mathcal P_{k,p,\alpha}$
		satisfies $d(G)\le (1-1/\rho)|V(G)|$. This holds when $G$ is edgeless.
		Deleting isolated vertices leaves $d(G)$ unchanged and decreases
		$|V(G)|$, so we may assume that $G$ has at least one edge and no
		isolated vertices. Write $d\coloneqq d(G)$ and $e\coloneqq |V(G)|-d$.
		The desired bound is equivalent to
		\begin{equation}\label{eq:pos-upper-target}
			|V(G)|\le \rho e.
		\end{equation}
		
		Fix a representation of $G$ as a subgraph of a host from
		Construction~\ref{cons:pos-amplification}. Write
		\[
		n_0\coloneqq |V|,\qquad q\coloneqq |Q|,\qquad
		u\coloneqq |U_0|,\qquad a\coloneqq\sum_{i=1}^p |U_i|.
		\]
		Then
		\begin{equation}\label{eq:pos-upper-order}
			|V(G)|=n_0+q+u+a,
			\qquad
			d=n_0+q+u+a-e.
		\end{equation}
		
		By Lemma~\ref{lem:pos-extension}, every independent set in $G$ has
		size at most $|V(G)|-d=e$. The sets $Q$ and $U_0$ are independent by
		Rules~\ref{item:pos-ampl-internal} and~\ref{item:pos-ampl-colour},
		respectively. When $k\ge 4$, the set $V$ is also independent, since
		each edge has at most three vertices in $V$ by
		Rule~\ref{item:pos-ampl-core}. Hence
		\begin{equation}\label{eq:pos-upper-parts}
			q\le e,\quad u\le e,\quad
			n_0\le e\quad\text{if }k\ge 4.
		\end{equation}
		Each colour block $U_{B_j}$ is independent by
		Rule~\ref{item:pos-ampl-block}, so $|U_{B_j}|\le e$ for every
		$j\in[p]$. By Lemma~\ref{lem:pos-blocks}, each colour $i\in[p]$
		belongs to exactly $m$ of the blocks $B_1,\ldots,B_p$. Thus the sum
		of the colour block sizes counts each vertex in $U_1\cup\cdots\cup U_p$
		exactly $m$ times, giving
		\begin{equation}\label{eq:pos-upper-blocks}
			ma=\sum_{j=1}^p |U_{B_j}|\le pe,
			\qquad\text{so}\qquad
			a\le \frac{p}{m}e.
		\end{equation}
		
		If $n_0+q+u\le 2e$, then \eqref{eq:pos-upper-order}
		and~\eqref{eq:pos-upper-blocks} give $|V(G)|\le (2+p/m)e$.
		This proves~\eqref{eq:pos-upper-target} in this case, since
		$\rho>2+p/m$. Indeed, by~\eqref{eq:pos-R}, the difference
		$\rho-(2+p/m)$ is $\alpha>0$ when $k\ge 4$, and is
		$\alpha/(1-\alpha)-1/m\ge 1/8-1/72>0$ when $k=3$.
		We may therefore assume that
		\begin{equation}\label{eq:pos-upper-large}
			n_0+q+u>2e.
		\end{equation}
		In particular, $n_0>0$, since $q,u\le e$ by~\eqref{eq:pos-upper-parts}.
		
		To bound $n_0+q$, let $H'$ be the $3$-graph on $V$ whose edges are
		the triples in $V$ covered in $G$. By Rule~\ref{item:pos-ampl-core},
		each such triple is an edge of $H$. Thus $H'\in\mathcal C_\theta$,
		since $\mathcal C_\theta$ is subgraph-closed. We verify that $H'$
		has no isolated vertices and that $d(H')>h(H')$, so that we can
		apply Lemma~\ref{lem:pos-core}.
		
		First, let $S\subseteq V$ be a pair covered in $G$. If there exists an edge
		containing $S$ and contained in $S\cup Q$, then it would lie in $K$ and
		have exactly two vertices in $V$, contrary to
		Rule~\ref{item:pos-ampl-internal}. This allows us to apply Lemma~\ref{lem:pos-extension} with $T=S$ and $I=Q$, and we obtain at least $d$ vertices
		$w\notin S\cup Q$ for which $S\cup\{w\}$ is covered in $G$.
		At most $a+u$ of these vertices lie in $U_0\cup\cdots\cup U_p$.
		Each remaining vertex lies in $V$ and gives an edge of $H'$ containing
		$S$ by the definition of $H'$. Hence
		\begin{equation}\label{eq:pos-upper-pair}
			|N_{H'}(S)|\ge d-a-u=n_0+q-e>e-u\ge 0,
		\end{equation}
		where the equality follows from~\eqref{eq:pos-upper-order}, the strict
		inequality from~\eqref{eq:pos-upper-large}, and the last inequality
		from~\eqref{eq:pos-upper-parts}. Additionally, every pair in $V$ covered in $G$ is also covered
		in $H'$. The converse follows from the definition of $H'$. Moreover, such a covered pair $S$ exists. Indeed, since $|V\cup U_0|=n_0+u>e$ by~\eqref{eq:pos-upper-large} and every independent set has size at most $e$ by Lemma~\ref{lem:pos-extension}, we derive that $V\cup U_0$ contains an edge of $G$. By
		Rule~\ref{item:pos-ampl-colour}, this edge has at most $k-2$ vertices
		in $U_0$, so it contains a pair in $V$. This fact, together with~\eqref{eq:pos-upper-pair}, implies
		\begin{equation}\label{eq:pos-lb-H'}
			d(H')\ge n_0+q-e>0.
		\end{equation}
		In particular, $H'$ is nonempty.
		
		Now fix $v\in V$. Since $G$ has no isolated vertices, $\{v\}$ is
		covered in $G$. If there exists an edge containing $v$ and contained in $\{v\}\cup U_0$, then it would have $k-1$ vertices of colour $0$, contrary to
		Rule~\ref{item:pos-ampl-colour}. Applying
		Lemma~\ref{lem:pos-extension} with $T=\{v\}$ and $I=U_0$, we obtain
		at least $d$ vertices $w\notin\{v\}\cup U_0$ for which $\{v,w\}$
		is covered in $G$. At most $a+q$ of them lie outside $V$, in
		$U_1\cup\cdots\cup U_p$ or in $Q$. Each remaining vertex lies in $V$ and is a neighbour
		of $v$ in $\partial_2H'$. It follows that
		\[
		\deg_{\partial_2H'}(v)\ge d-a-q=n_0+u-e.
		\]
		This is positive, since~\eqref{eq:pos-upper-large} and~\eqref{eq:pos-upper-parts} give $n_0+u>2e-q\ge e$. Thus $H'$ has no isolated
		vertices. We also obtain
		\begin{equation}\label{eq:pos-upper-shadow}
			h(H')=n_0-\delta(\partial_2H')
			\le n_0-(n_0+u-e)=e-u.
		\end{equation}
		
		Using~\eqref{eq:pos-lb-H'} and~\eqref{eq:pos-upper-shadow}, we obtain
		\[
		d(H')\ge n_0+q-e>e-u\ge h(H'),
		\]
		where the strict inequality follows from~\eqref{eq:pos-upper-large}.
		Lemma~\ref{lem:pos-core} now gives
		$d(H')\le n_0/(9-\theta)=\alpha n_0$. Combining the lower and upper
		bounds on $d(H')$, we obtain
		\begin{equation}\label{eq:pos-upper-core-size}
			n_0+q-e\le\alpha n_0,
			\quad\text{or equivalently}\quad
			n_0+q\le e+\alpha n_0.
		\end{equation}
		
		Suppose first that $k\ge 4$. By~\eqref{eq:pos-upper-parts},
		we have $n_0\le e$, so~\eqref{eq:pos-upper-core-size} gives
		\[
		n_0+q\le e+\alpha n_0\le (1+\alpha)e.
		\]
		Adding $u\le e$ from~\eqref{eq:pos-upper-parts} and
		$a\le (p/m)e$ from~\eqref{eq:pos-upper-blocks}, we obtain
		\[
		|V(G)|=(n_0+q)+u+a
		\le (1+\alpha)e+e+\frac{p}{m}e
		=\rho e.
		\]
		
		Finally, suppose that $k=3$. Then $q=0$ by construction, so
		\eqref{eq:pos-upper-core-size} gives $n_0\le e/(1-\alpha)$.
		For each $i\in[p]$, every triple in $U_0\cup U_i$ contains two
		vertices of the same colour and is therefore forbidden by
		Rule~\ref{item:pos-ampl-colour}. Thus $U_0\cup U_i$ is independent,
		and Lemma~\ref{lem:pos-extension} gives $u+|U_i|\le e$.
		Summing up over $i\in [p]$ and dividing by $p$ yields $u+a/p\le e$. Using this inequality
		and then~\eqref{eq:pos-upper-blocks}, we obtain
		\[
		u+a
		=\left(u+\frac ap\right)+\left(1-\frac1p\right)a
		\le e+\left(1-\frac1p\right)\frac{p}{m}e
		=\left(1+\frac{p-1}{m}\right)e.
		\]
		Adding the bound for $n_0$ gives
		\[
		|V(G)|=n_0+u+a
		\le \left(\frac{1}{1-\alpha}+1+\frac{p-1}{m}\right)e
		=\rho e.
		\]
		This proves~\eqref{eq:pos-upper-target} and completes the upper bound.
		
		\noindent\textbf{The lower bound.}
		Fix a sufficiently large integer $n$ and put $t\coloneqq n/\rho$.
		For $1\le i\le p$, take $|U_i|=\lfloor t/m\rfloor$. Choose
		$q=|Q|$ and $u=|U_0|$ as in the table below, and put the remaining
		$n_0=n-q-u-p\lfloor t/m\rfloor$ vertices in $V$. The last column
		follows from the definition of $\rho$ in~\eqref{eq:pos-R}.
		\begin{equation}\label{eq:pos-lower-sizes}
			\begin{array}{c|ccc}
				& q & u & n_0 \\ \hline
				k=3    & 0 & \lfloor(1-1/m)t\rfloor & t/(1-\alpha)+O(1) \\
				k\ge 4 & \lfloor\alpha t\rfloor & \lfloor t\rfloor & t+O(1)
			\end{array}
		\end{equation}
		Here $k,p$ and $\alpha$ are fixed and $n\to\infty$. In particular,
		every outside colour class $U_i$ has size at most $t$, and each colour block has
		size
		\begin{equation}\label{eq:pos-lower-block-size}
			|U_{B_j}|=m\lfloor t/m\rfloor\le t
			\qquad (j\in[p]).
		\end{equation}
		
		By Lemma~\ref{lem:pos-core}, we may choose $H\in\mathcal C_\theta$
		on $V$ with minimum codegree
		$\delta_2(H)\ge\alpha n_0-2$. Let $G$ be the corresponding host
		from Construction~\ref{cons:pos-amplification}, with every permitted
		edge included. For sufficiently large $n$, the graph $H$ is nonempty
		and $q\ge k-3$. Any edge of $H$ together with $k-3$ vertices of $Q$
		forms an edge of $G$, so $G$ is nonempty.
		
		Fix a facet $F$ of $G$, namely $F$ is a covered $(k-1)$-set of $G$. We will prove that the number of non-neighbours of $F$ is at most
		$t+O(1)$. We split the proof into three cases according to $|F\cap K|$.
		
		\smallskip
		\noindent\textbf{Case 1.} $|F\cap K|\ge 2$.
		The facet $F$ contains at most $k-3$ outside vertices. For any
		outside vertex $x\notin F$, the set $F\cup\{x\}$ meets both $K$ and
		its complement and contains at most $k-2$ outside vertices in total.
		Thus every
		outside vertex $x\notin F$ is a neighbour of $F$.
		
		We next consider vertices in $K$. Suppose first that $F\cap V$ is
		a pair $S$. Then $N_H(S)\subseteq N_G(F)$, since for each
		$x\in N_H(S)$ the three vertices of $F\cup\{x\}$ in $V$ form an
		edge of $H$. Hence every non-neighbour of $F$ outside $F$ lies in
		$Q\cup(V\setminus N_H(S))$, whose size is at most
		\[
		q+n_0-\delta_2(H)
		\le q+(1-\alpha)n_0+2
		=t+O(1).
		\]
		The last equality follows from either row of~\eqref{eq:pos-lower-sizes}.
		
		When $k=3$, this covers all possibilities in Case~1, since $K=V$
		and $F$ is a pair. Now suppose that $k\ge 4$. The remaining
		possibilities for $|F\cap V|$ are $0$, $1$ and $3$. If $F\cap V$
		is empty, then every vertex of $V$ is a neighbour of $F$, since
		adding it gives an edge with exactly one vertex in $V$. Thus only
		vertices of $Q$ can be non-neighbours outside $F$, and there are at most $q$
		of them. If $|F\cap V|$ is one or three, then every vertex of
		$Q\setminus F$ is a neighbour of $F$, since adding it preserves
		$F\cap V$ and its size is allowed by Rule~\ref{item:pos-ampl-internal}.
		In this case only vertices of $V$ can be non-neighbours outside $F$, and there
		are at most $n_0$ of them. Both bounds are at most $t+O(1)$
		by~\eqref{eq:pos-lower-sizes}.
		
		\smallskip
		\noindent\textbf{Case 2.} $|F\cap K|=1$.
		Here $F$ has $k-2$ outside vertices. Every $x\in K\setminus F$ is
		a neighbour of $F$, since $F\cup\{x\}$ meets both $K$ and its
		complement, has at most two vertices in $V$, and has the same outside
		vertices as $F$.
		
		For an outside vertex $x\notin F$, the only rule that can exclude
		$F\cup\{x\}$ is Rule~\ref{item:pos-ampl-colour}. This happens
		precisely when all $k-2$ outside vertices of $F$ lie in a single
		class $U_i$ and $x\in U_i$. Thus all non-neighbours outside $F$, if any, lie in one outside colour class, and its size is at most $t$.
		
		\smallskip
		\noindent\textbf{Case 3.} $F\cap K=\varnothing$.
		Let $J$ be the set of colours appearing in $F$. Since $F$ is covered, it lies in some edge $E$ of $G$.
		By Rule~\ref{item:pos-ampl-colour}, each outside colour class contains
		at most $k-2$ vertices of $E$. Thus $|J|\ge 2$. 
		In this case, every vertex of $K$ is a neighbour of $F$, since adding it gives a
		$k$-set meeting $K$ in exactly one vertex and preserves the outside vertices.
		
		Suppose first that $J\subseteq B_j$ for some block $B_j$. This block
		is unique, as any two distinct colours in $J$ belong together to
		at most one block by Lemma~\ref{lem:pos-blocks}. An outside vertex
		$x\notin U_{B_j}$ has a colour outside of $J$, so
		$F\cup\{x\}$ satisfies Rule~\ref{item:pos-ampl-colour}. It also
		satisfies Rule~\ref{item:pos-ampl-block}, since $B_j$ is the only
		block containing $J$.
		Thus $x$ is a neighbour of $F$. All non-neighbours of $F$ therefore lie in
		$U_{B_j}$, which has size at most $t$
		by~\eqref{eq:pos-lower-block-size}.
		
		Finally, suppose that no block contains $J$. Then $F\cup\{x\}$
		satisfies Rule~\ref{item:pos-ampl-block} for every outside vertex
		$x\notin F$. Such a vertex can be a non-neighbour of $F$ only if
		its colour already appears $k-2$ times in $F$. For $k\ge 4$, there
		is at most one such colour and its class
		has size at most $t$. For $k=3$, the two vertices of $F$ have
		distinct colours, and all non-neighbours lie in the union of their
		two classes. If one of the colours is $0$, the two classes have
		total size at most $(1-1/m)t+t/m=t$. If both colours are nonzero,
		their total size is at most $2t/m\le t$.
		
		In all three cases, at most $t+O(1)$ vertices of $V(G)\setminus F$
		lie outside $N_G(F)$. Since $|F|=k-1$ is fixed and $t=n/\rho$, taking the minimum over all facets gives
		\[
		d(G)\ge n-(k-1)-t-O(1)
		=\left(1-\frac{1}{\rho}\right)n-O(1).
		\]
		Since the construction works for every sufficiently large
		$n$, this lower bound and the upper bound together show that
		$\mathcal P_{k,p,\alpha}$ has positive codegree density $1-1/\rho$.\end{proof}
		
		\section*{AI Disclosure} AI tools (ChatGPT 6.0 Astra) generated the main proof ideas in this project. Through repeated discussions with the AI, we obtained simpler proofs. We then revised the proofs based on our understanding and simplified them further. AI also assisted with writing and revising the manuscript and with drawing Figure~\ref{fi:Consruction}. We have independently checked all proofs and take full responsibility for the content of this paper.
			
	\bibliographystyle{abbrv}%abbrv
	\bibliography{codegree}

\begin{thebibliography}{10}

\bibitem{BT11}
R.~Baber and J.~Talbot.
\newblock Hypergraphs do jump.
\newblock {\em Combin. Probab. Comput.}, 20(2):161--171, 2011.

\bibitem{BaloghHalfpapLidickyPalmer26}
J.~Balogh, A.~Halfpap, B.~Lidick{\'y}, and C.~Palmer.
\newblock Positive co-degree densities and jumps.
\newblock {\em Innov. Graph Theory}, 3:1--36, 2026.

\bibitem{E64}
P.~Erd{\H{o}}s.
\newblock On extremal problems of graphs and generalized graphs.
\newblock {\em Israel J. Math.}, 2:183--190, 1964.

\bibitem{ES66}
P.~Erd{\H {o}}s and M.~Simonovits.
\newblock A limit theorem in graph theory.
\newblock {\em Studia Sci. Math. Hungar.}, 1:51--57, 1966.

\bibitem{ES46}
P.~Erd{\H {o}}s and A.~H. Stone.
\newblock On the structure of linear graphs.
\newblock {\em Bull. Amer. Math. Soc.}, 52:1087--1091, 1946.

\bibitem{FPRT07}
P.~Frankl, Y.~Peng, V.~R\"odl, and J.~Talbot.
\newblock A note on the jumping constant conjecture of {E}rd{\H {o}}s.
\newblock {\em J. Combin. Theory Ser. B}, 97(2):204--216, 2007.

\bibitem{FR84}
P.~Frankl and V.~R\"odl.
\newblock Hypergraphs do not jump.
\newblock {\em Combinatorica}, 4(2-3):149--159, 1984.

\bibitem{GPRS26}
J.~Gao, O.~Pikhurko, M.~Rong, and S.~Sun.
\newblock Rational codegree {Tur{\'a}n} density of hypergraphs.
\newblock Preprint, arXiv:2601.00758, 2026.

\bibitem{HLP25}
A.~Halfpap, N.~Lemons, and C.~Palmer.
\newblock Positive co-degree density of hypergraphs.
\newblock {\em J. Graph Theory}, 110(2):209--222, 2025.

\bibitem{JLR00}
S.~Janson, T.~{\L}uczak, and A.~Ruci{\'n}ski.
\newblock {\em Random Graphs}.
\newblock Wiley-Interscience Series in Discrete Mathematics and Optimization. John Wiley \& Sons, 2000.

\bibitem{kee11}
P.~Keevash.
\newblock Hypergraph {T}ur\'an problems.
\newblock In {\em Surveys in combinatorics 2011}, volume 392 of {\em London Math. Soc. Lecture Note Ser.}, pages 83--139. Cambridge Univ. Press, Cambridge, 2011.

\bibitem{LiLiu26v2}
H.~Li and X.~Liu.
\newblock On possible {Tur{\'a}n} densities under general uniformity conditions.
\newblock Preprint, arXiv:2610.07049v2, 2026.

\bibitem{LiLiu2026}
H.~Li and X.~Liu.
\newblock On the spectrum of minimum {$\ell$}-degree {Tur{\'a}n} densities.
\newblock Preliminary manuscript, August 26, 2026.

\bibitem{LinZhou2026}
H.~Lin and W.~Zhou.
\newblock The {$\ell$}-degree {Tur{\'a}n} densities of hypergraphs: infinite families and their finite subfamilies.
\newblock Preprint, arXiv:2610.17558, 2026.

\bibitem{LM26}
X.~Liu and D.~Mubayi.
\newblock The number {$4/9$} is a non-jump for {$3$}-graphs.
\newblock Preprint, {arXiv}:2605.13567, 2026.

\bibitem{LP26}
X.~Liu and O.~Pikhurko.
\newblock Intervals of hypergraph {Tur{\'a}n} densities.
\newblock Preprint, {arXiv}:2605.25914, 2026.

\bibitem{MZ07}
D.~Mubayi and Y.~Zhao.
\newblock Co-degree density of hypergraphs.
\newblock {\em J. Combin. Theory Ser. A}, 114(6):1118--1132, 2007.

\bibitem{P14}
O.~Pikhurko.
\newblock On possible {T}ur\'an densities.
\newblock {\em Israel J. Math.}, 201(1):415--454, 2014.

\bibitem{Pikhurko23}
O.~Pikhurko.
\newblock On the limit of the positive {$\ell$}-degree {T}ur\'an problem.
\newblock {\em Electron. J. Combin.}, 30(3):Paper No. 3.25, 15, 2023.

\bibitem{Shaw26}
B.~Randall~Shaw.
\newblock Minimal hypergraph non-jumps.
\newblock {\em European J. Combin.}, 137:12, 2026.
\newblock Id/No 104400.

\bibitem{sido95}
A.~Sidorenko.
\newblock What we know and what we do not know about {T}ur\'an numbers.
\newblock {\em Graphs Combin.}, 11(2):179--199, 1995.

\bibitem{YP23}
Z.~Yan and Y.~Peng.
\newblock Non-jumping {Tur{\'a}n} densities of hypergraphs.
\newblock {\em Discrete Math.}, 346(1):11, 2023.
\newblock Id/No 113195.

\end{thebibliography}
\end{document}